\documentclass[twoside,6pt,leqno]{article}
\usepackage{textcase}
\usepackage{amsmath}
\usepackage{amssymb}
\usepackage{amsbsy}
\usepackage{amscd}
\usepackage{amsfonts}
\usepackage{amsopn}
\usepackage{amstext}
\usepackage{amsthm}
\usepackage{amsxtra}
\usepackage{graphicx,setspace}

\newtheorem{theorem}{Theorem}[section]
\newtheorem{lemma}{Lemma}[section]
\newtheorem{proposition}{Proposition}[section]
\newtheorem{corollary}{Corollary}[section]

\def\reverseddots{\mathinner{\mkernlmu\raise7pt
\vbox{\kern7pt\hbox{.}}
\mkern-12mu\raise4pt\hbox{.}
\mkern-12mu\raise1pt\hbox{.}
\mkern1mu}}
\title{\bf{The central value of the Rankin-Selberg $L$-functions}}
\author{Xiaoqing Li}
\date{}
\begin{document}
\maketitle
\pagenumbering{arabic}
\pagestyle{myheadings}
\markboth{}{}
\thispagestyle{headings}
\begin{abstract} Let $f$ be a Maass form for $SL(3, \mathbb{Z})$ which is fixed and $u_j$ be an orthonormal basis of even Maass
forms for $SL(2, \mathbb{Z}),$ we prove an asymptotic formula for the average of the product of the Rankin-Selberg $L$-function of 
$f$ and $u_j$ and the $L$-function of $u_j$ at the central value $1/2$. This implies simultaneous nonvanishing 
results of these $L$-functions at $1/2.$ 
\end{abstract}
%\begin{center}
%\center{
\section
{Introduction}
The values of $L$-functions at special points have been the subject of intensive studies. For example, a good positive lower bound
for the central value of Hecke $L$-functions would rule out the existence of the Landau-Siegel zero, see the notable paper [IS]; the 
nonvanishing of certain Rankin-Selberg $L$-functions is a crucial ingredient in the current development of the generalized Ramanujan conjecture [LRS], 
etc. In this paper, we consider the simultaneous nonvanishing problem of products of Rankin-Selberg on ${\rm GL}(3)$ and 
${\rm GL}(2)$ 
and Maass $L$-functions on ${\rm GL}(2)$ at the central point $1/2$. \\
Specifically, let $u_j(z)$ be an orthonormal basis of even Hecke-Maass forms for the modualr group $SL(2, \mathbb{Z}).$ 
For each $u_j(z),$ let $a_j(n)$ be its normalized Fourier coefficients (see the next section), we associate the $L$-function:
\begin{equation}
L(s, u_j)=\sum\limits_{n\geqslant 1}a_j(n)n^{-s}
\end{equation}
which has analytic continuation to the whole complex plane and  satisfies a functional equation relating $s$ to $1-s.$
Let $f(z)$ be a Hecke-Maass form of type $(\nu_1, \nu_2)$ for $SL(3, \mathbb{Z})$ and $\tilde{f}(z)$ be its dual Maass form. $f(z)$
has a Fourier-Whittaker expansion with Fourier coefficients $A(m, n).$ The $L$-function
\begin{equation}
L(s, f)=\sum\limits_{m=1}^{\infty} A(1, m)m^{-s}
\end{equation}
has analytic continuation to the whole complex plane and satisfies a functional equation. The Rankin-Selberg $L$-function defined by
\begin{equation}
L(s, f\times u_j)=
\sum\limits_{m\geqslant 1}\sum\limits_{n\geqslant 1}\frac{\bar {a}_j(n)A(m, n)}{(m^2n)^s}
\end{equation}
also has analytic continuation to the whole complex plane and satisfies a functional equation relating $s$ to $1-s.$ See the next section for related terminology and details.\\
Our main theorem is the following:
\begin{theorem}
For $f$ a fixed Hecke-Maass form for $SL(3, \mathbb{Z})$ and $\tilde{f}$ be its dual Maass form, $u_j$ an orthonormal basis
 of even Hecke-Maass forms of type $\frac{1}{2}+it_j$ for $SL(2, \mathbb{Z}),$
we have
\begin{equation}
{\sum_j}^{'} e^{-\frac{t_j^2}{T^2}}L\left(\frac{1}{2}, f\times u_j\right)L\left(\frac{1}{2}, u_j\right)
=\frac{12 L(1, f)L(1, \tilde{f})}{\pi^3}T^2
+O_{\varepsilon, f}(T^{\frac{11}{6}+\varepsilon})
\end{equation}
where ' means summing over the orthonormal basis of even Hecke-Maass forms and $\varepsilon>0$ is arbitrarily small.
\end{theorem}
It is known [JS] that $L(1, f)L(1, \tilde{f})\neq 0,$ so we have
\begin{corollary}
Under the same assumption as in the above theorem, there are infinitely many $u_j^{'}$s such that
$$
L\left(\frac{1}{2}, f\times u_j\right)L\left(\frac{1}{2}, u_j\right)\neq 0.
$$
\end{corollary}
\noindent{\bf Remarks} 1. If $f$ comes from the Gelbart-Jacquet lift [GJ] from ${\rm GL}(2),$ then there is Watson's formula [Wa] which relates
$L\left(\frac{1}{2}, f\times u_j\right)L\left(\frac{1}{2}, u_j\right)$ to some period integrals. Then the nonvanishing of such $L$-functions
at the central point implies the nonvanishing of those periods, see also [Re] and [GJR] from the representation theory point
 of view.\\
2. The technology in this paper also yields 
$$
{\sum_j}^{'} e^{-\frac{(t_j-T)^2}{H^2}}L\left(\frac{1}{2}, f\times u_j\right)L\left(\frac{1}{2}, u_j\right)
= O_{\varepsilon, f}(T^{\frac{11}{6}+\varepsilon})
$$
where $H=T^{\frac{5}{6}},$ ' means summing over the orthonormal basis of even Hecke-Maass forms. When $f$ is selfdual, 
by the positivity of the $L$-functions
([La], [KS], [Gu]), we have
$$
L\left(\frac{1}{2}, f\times u_j\right)L\left(\frac{1}{2}, u_j\right)
\ll_{\varepsilon, f}T^{\frac{11}{6}+\varepsilon}
$$
for $t_j-T\asymp H.$ This yields the subconvexity of the product of the $L$-functions which is as strong as the current record 
subconvexity bound  
\begin{equation}
L\left(\frac{1}{2}, u_j\right)\ll \left(1+|t_j|\right)^{\frac{1}{3}}
\end{equation}
combining with the convexity bound 
\begin{equation}
L\left(\frac{1}{2}, f\times u_j\right)\ll \left(1+|t_j|\right)^{\frac{3}{2}}.
\end{equation}
(1.5) was first proved conditionally by Iwaniec in [Iw1] and an unconditional proof was given 
by Ivic [Iv1]
and subsequently by Jutila [Ju], while the convexity bound (1.6) remains untouched. In the case that 
$f$ comes from the Gelbart-Jacquet lift [GJ] from ${\rm GL}(2),$ Bernstein
and Reznikov [BR] obtained the bound $L\left(1/2, f\times u_j\right)
L\left(1/2, u_j\right)\ll
\left(1+|t_j|\right)^{\frac{5}{3}+ \varepsilon}$ using the representation theory for compact Riemann sufaces and they claimed their method should also work 
in general.\\
3. Much stronger nonvanishing results in terms of percentage of nonvanishing are known for lower degree $L$-functions using
 the powerful mollification techniques, see [IS], [Lu], [KMV], [So], for example. In our case, such stronger results haven't been
done yet.\\
Our approach to prove Theorem 1.1 makes use of the Kuznetsov formula on ${\rm GL}(2)$ and the Voronoi formula on ${\rm GL}(3)$ which was first
derived by Miller and Schmidt [MS1] using the theory of automorphic distributions, see also [GL] for a simple, analytic proof. 
The Voronoi formula on ${\rm GL}(3)$ has been used by Sarnak and Watson, Miller and Schmidt (see [Mi], [MS2]) to prove a variety of results on 
$L$-functions, our paper gives another application of this very useful tool.
 \section {A review of automorphic forms}
\setcounter{equation}{0} 
We set up the problem in a general background.\\
For $n\geqslant 2,$ let $G={\rm GL}(n, \mathbb{R}), \Gamma=SL(n, \mathbb{Z})$ and 
$$
\frak h^n={\rm GL}(n, \mathbb{R})/\langle O(n, \mathbb R)\cdotp\mathbb {R}^{\times}\rangle
$$
be the generalized upper half plane. Every element $z\in \frak{h^n}$ has the form $z=xy$ where
$$
x=
\left(\!
\begin{matrix}
1&x_{1,2}&x_{1,3}&\ldots&x_{1,n}\\
&1&x_{2,3}&\ldots&x_{2,n}\\
&&\ddots&&\vdots\\
&&&1&x_{n-1,n}\\
&&&&1
\end{matrix}
\!\!\right),
$$
%\end{document}
$$
y=\mbox{diag}(y_1y_2...y_{n-1},\; y_1y_2...y_{n-2},\; ..., y_1,\;
1),
$$
with $x_{ij}\in \mathbb {R}$ for $1\leqslant i< j\leqslant n$ and $y_i>0$ for $1\leqslant i\leqslant n-1.$

Let $\nu=(\nu_1, \nu_2, ..., \nu_{n-1})\in \mathbb {C}^{n-1}.$ The function
$$
I_{\nu}(z)=\prod_{i=1}^{n-1}\prod_{j=1}^{n-1}y_i^{b_{n-i, j}\nu_j}
$$
with
$$
b_{i, j}=
\left\{
\begin{array}{lll}
ij& \mbox{} &\text{if}\; i+j\leqslant n, \\
(n-i)(n-j)&\mbox{}&  \text{otherwise,} \end{array} \right.
$$
is an eigenfunction of every differential operator $D$ in ${\cal {D}}^n,$ the center of the universal enveloping algebra of 
$gl(n, \mathbb{R}).$ Here $gl(n, \mathbb{R})$ is the Lie algebra of ${\rm GL}(n, \mathbb{R}).$ Let us write
$$
DI_{\nu}(z)=\lambda_D I_{\nu}(z)
$$
for every $D\in {\cal {D}}^n.$ An automorphic form $f$ of type $\nu$ for $\Gamma=SL(n, \mathbb{Z})$ is a smooth
 function on $\frak{h^n}$which satisfies\\

1)\; $f(\gamma z)=f(z)$ for all $\gamma\in \Gamma;$\\

2)\; $Df(z)=\lambda_D f(z)$ for all $D\in {\cal {D}}^n.$\\

\noindent If $f$ also satisfies\\

3)\; $\int\limits_{\Gamma\cap U \backslash U}f(uz)d^{*}u=0$\\
 
\noindent where $d^{*}u=\prod\limits_{1\leqslant i<j\leqslant n}du_{i, j},$
$U$ is formed by all upper triangular matrices of the form
$$u=\left(\!
\begin{matrix}
I_{r_1}&&&&\\
&I_{r_2}&&*&\\
&&\ddots&&\\
&&&&I_{r_m}
\end{matrix}
\!\right),
$$
with $r_1+r_2+\cdots+r_m=n,$ $I_r$ denotes the $r\times r$ identity matrix and $*$ denotes arbitrary real elements, then $f$ is called a Maass form of type $\nu.$ 

For $z\in \frak{h^n},$ let $U_n({\mathbb R})$ denote the group of $n\times n$ upper triangular matrices
with ones on the diagonal. Let
$$
W_{\rm Jacquet}(z; \nu, \psi_m)=\int\limits_{U_n(\mathbb{R})}I_{\nu}(w_nuz)\overline{\psi_m(u)}
d^{*}u
$$ 
be Jacquet's 
Whittaker function which has rapid decay as 
$y_i\rightarrow \infty, 1\leqslant i\leqslant n-1.$
Here
$$
\psi_m(u)=e(m_1u_{1,2}+m_2u_{2, 3}+\cdots+m_{n-1}u_{n-1, n})
$$ with $e(z):=e^{2\pi iz}$ throughout the paper 
and
$$w_n=\left(\!
\begin{matrix}
&&&&&\pm 1\\
&&&&1&\\
%&&\reverseddots&&&\\
&& ...&&&\\

1&&&&&
\end{matrix}
\!\right).
$$
Every Maass form $f(z)$ of type $\nu=(\nu_1, ...,
\nu_{n-1})$ has the following Fourier-Whittaker expansion:
\begin{eqnarray}
\lefteqn{\;\;\;\;\;\;\;\;\;\;\;\;\;\;
f(z)=
%\label{(2.4)}}
\!\!\!\sum\limits_{\gamma\in U_{n-1}(\mathbb{Z})\backslash SL(n-1, \mathbb{Z})}
\;\sum\limits_{m_1=1}^{\infty}\cdots\sum\limits_{m_{n-2=1}}^{\infty}
\sum\limits_{m_{n-1}\neq 0}
\frac{A(m_1, \ldots, m_{n-1})}{\prod\limits_{k=1}^{n-1}|m_k|^{\frac{k(n-k)}{2}}}
\;\;\;\;\;\;\;\;\;\;\;\;\;\;\label{(2.5)}}\\
&&\hspace{4cm}\;\;\;\;\;\;\cdot W_{\rm Jacquet}\left(M\left(\!\!\begin{array}{rr}
                   \gamma\!&\! \\
                   \!&\! 1\\
                   \end{array}\!\!\right)z, \nu, \psi_{1,\cdots, 1, \frac{m_{n-1}}{|m_{n-1}|}}\right)\nonumber
,\end{eqnarray}
where $U_n(\mathbb Z)$ is the subgroup of $U_n(\mathbb R)$ with coefficients in 
$\mathbb Z$, and \\
$M=\mbox{diag}\left(m_1\cdots m_{n-2}|m_{n-1}|,\cdots, m_1m_2, m_1, 1\right).$ It is easy to prove
that (see Chapter 9 in [Go]) the
dual Maass form $\tilde{f}(z):=f(w_n {^tz^{-1}}w_n)$ is a Maass form of 
type $(\nu_{n-1},\cdots
,\nu_1)$ with Fourier coefficients $A(m_{n-1}, \ldots, m_1).$ 

Next let's recall some facts about Hecke operators. Let
 ${\cal L}^2(\Gamma\setminus \frak{h^n})$ be the space 
of square integrable
automorphic forms for $\Gamma$ equipped with the inner product:
$$
\langle f, g\rangle=\int\limits_{\Gamma\setminus \frak{h^n}}f(z)\overline{g(z)}\;d^*(z),$$
for all $f, g\in {\cal L}^2(\Gamma\setminus \frak{h^n}),$ where $d^*(z)=\prod\limits_{1\leqslant i<j\leqslant n}
dx_{i, j}\prod\limits_{k=1}^{n-1}
y_k^{-k(n-k)-1}dy_k$ is the $G$ left invariant measure. For every integer $N\geqslant 1,$ we define a Hecke operator
$T_N$ acting on ${\cal L}^2(\Gamma\setminus \frak{h^n})$ by the following formula:
%\end{document}
$$
T_Nf(z)=\frac{1}{N^{\frac{n-1}{2}}}\!\!\!\!
\sum_{\substack{
\prod\limits_{l=1}^{n}
c_l=N\\0\leqslant c_{i, l}<c_l\; 
(1\leqslant i <l\leqslant n)}}f\left(
\left(\!
\begin{matrix}
c_1&c_{1,2}&\ldots&c_{1,n}\\
&c_2&\ldots&c_{2,n}\\
&&\ddots&\vdots\\
&&&c_n
\end{matrix}
\!\!\right)\cdotp z
\right).
$$ 
%\end{document}
The Hecke operators are normal operators. They commute with each other as well as with the $G$ 
invariant differential operators.  So we may simultaneously diagonalize the space ${\cal L}^2
(\Gamma\setminus \frak{h^n})$
by all these operators. Let $f$ be a Maass form with Fourier expansion (2.1) which is also 
an eigenfunction of all the 
Hecke operators. We normalize $A(1, \ldots, 1)$ to be 1. Then we have the following multiplicativity 
relations:
$$
A(m_1m_1^{\prime},\;\ldots,\; m_{n-1}m_{n-1}^{\prime})=
A(m_1,\; \ldots,\; m_{n-1})\cdot A(m_1^{\prime},\;\ldots, \;m_{n-1}^{\prime}),
$$
%\end{document}
if $(m_1\ldots m_{n-1}, m_1^{\prime}\ldots m_{n-1}^{\prime})=1,$ and 
$$
A(m, 1, \ldots,  1)A(m_1,\ldots, m_{n-1})=\!\!\!\!\!\!\!\!\!\!\!\!\!\!\!\!\!\!\!\!\!\!
\sum_{\substack{\prod\limits_{l=1}^{n}c_l=m\\
c_1|m_1,\; c_2|m_2,\; \ldots, \;c_{n-1}|m_{n-1}}}\!\!\!\!\!\!\!\!\!\!\!\!\!\!\!\!\!\!\!\!\!
A\left(\frac{m_1c_n}{c_1}, \frac{m_2c_1}{c_2},\ldots, \frac{m_{n-1}c_{n-2}}
{c_{n-1}}\right).
$$
The above material is taken from [Go]. Our main interests in this paper are the cases when $n=2$ and $3.$\\
For $n=2,$ one can identify $\frak{h^2}$ with the upper half plane $\mathbb {H}.$ ${\cal{D}}^2$ is generated by the Laplace operator
$$
\Delta=-y^2\left(
\frac{\partial^2}{\partial x^2}+\frac{\partial^2}{\partial y^2}\right)
$$
which has a spectral decomposition on $L^2(SL(2, \mathbb{Z})\setminus\mathbb{H}):$
$$
L^2(SL(2, \mathbb{Z})\setminus\mathbb{H})={\cal{C}}\oplus C(SL(2, \mathbb{Z})\setminus{\mathbb{H}})\oplus
{\cal{E}}(SL(2, \mathbb{Z})\setminus\mathbb{H}).
$$
Here $\cal {C}$ is the space of constant functions. $C(SL(2, \mathbb{Z})\setminus{\mathbb{H}})$ is the space of Maass forms 
and ${\cal{E}}(SL(2, \mathbb{Z})\setminus\mathbb{H})
$ is the space of Eisenstein series.\\
Let ${\cal{U}}=\{{u_j: j\geqslant 1}\}$ be an orthonormal basis of Hecke-Maass forms of type $s_j=\frac{1}{2}+it_j$ with $t_j\geqslant 0$in the 
space
$C(SL(2, \mathbb{Z})\setminus\mathbb{H}).$ Any $u_j(z)$ has the Fourier expansion
\begin{equation}
u_j(z)=\sum\limits_{n\neq 0}\rho_j(n)W_{s_j}(nz)
\end{equation}
where $W_s(z)$ is the Whittaker function given by
$$
W_s(z)=2|y|^{\frac{1}{2}}K_{s-\frac{1}{2}}(2\pi|y|)e(x)
$$
and $K_s(y)$ is the $K$-Bessel function. $C(SL(2, \mathbb{Z}\setminus\mathbb{H})$ consists of even Maass forms and odd Maass forms
according to $u_j(-\bar{z})=u_j(z)$ or $u_j(-\bar{z})=-u_j(z).$ The Eisenstein series $E(z, s)$ defined by
\begin{equation}
E(z, s)=\frac{1}{2}\sum_{\substack{c, d\in \mathbb{Z}\\ (c, d)=1}}\frac{y^s}{|cz+d|^{2s}}
\end{equation}
has the following Fourier expansion
$$
E(z, s)=y^s+\phi(s)y^{1-s}+\sum\limits_{n\neq 0}\phi(n, s)W_s(nz)
$$
where
$$
\phi(s)=\sqrt{\pi}\frac{\Gamma(s-\frac{1}{2})}{\Gamma(s)}\frac{\zeta(2s-1)}{\zeta(2s)}
$$ with $\zeta(s)$ be the Riemann zeta function
and
$$
\phi(n, s)=\pi^s\Gamma(s)^{-1}\zeta(2s)^{-1}|n|^{-\frac{1}{2}}\eta(n, s)
$$
with
$$
\eta(n, s)=\sum\limits_{ad=|n|}\left(\frac{a}{d}\right)^{s-\frac{1}{2}}.
$$
Let
\begin{equation}
a_j(n)=2\rho_j(n)|n|^{\frac{1}{2}}\Gamma\left(\frac{1}{2}+it_j\right),
\end{equation} we have Kuznetsov's formula (see [CI])
\begin{eqnarray}
\lefteqn{\;\;\;\;\;\;\;\;\;\;\;{\sum\limits_{j\geqslant 1}}^{'}h(t_j)\bar{a}_j(n)a_j(l)+\frac{1}{4\pi}\int\limits_{-\infty}^{\infty}h(r)\omega(r)
\bar{\eta}\left(n, \frac{1}{2}+ir\right)\eta\left(l, \frac{1}{2}+ir\right)dr}\\
&\;\;\;\;\;\;\;=\frac{1}{2}\delta(n, l)H
+\sum\limits_{c>0}\frac{1}{2c}\left\{S(n, l; c)H^+\left(\frac{2\sqrt{nl}}{c}\right)+
 S(-n, l; c)H^-\left(\frac{2\sqrt{nl}}{c}\right)\right\},\nonumber
\end{eqnarray}
where ${\sum}^{'}$ restricts to the even Maass forms, $\delta(n, l)$ is the Kronecker symbol,
\begin{equation}
H=\frac{1}{\pi}\int\limits_{-\infty}^{\infty}h(t)\tanh (\pi t) tdt,
\end{equation}
\begin{equation}
\omega(r)=4\pi\Bigg|\phi\left(1, \frac{1}{2}+ir\right)\Bigg|^2\cosh^{-1}\pi r,
\end{equation}
\begin{equation}
H^+(x)=2i\int\limits_{-\infty}^{\infty}J_{2it}(2\pi x)\frac{h(t)t}{\cosh(\pi t)}dt,
\end{equation}
\begin{equation}
H^-(x)=\frac{4}{\pi}\int\limits_{-\infty}^{\infty}K_{2it}(2\pi x)\sinh ({\pi}t) h(t)tdt,
\end{equation}
\begin{equation}
S(n, l; c)=\displaystyle{\sum_{d\bar {d}\equiv1(\text{mod}\; c)}}e\Big(\frac{dl+\bar{d}n}{c}\Big)
\end{equation}
is the classical Kloosterman sum.
(2.5) holds for any $n, l\geqslant 1$ and any test function $h(t)$ which is even and satisfies the following conditions:
\\
I) $h(t)$ is holomorphic in $|\Im t|\leqslant \sigma;$\\
II)$h(t)\ll (|t|+1)^{-\theta}$
for some $\sigma>\frac{1}{2}$ and $\theta >2.$\\
Now for $(\nu_1, \nu_2)\in {\cal C}^2,$ set 
\begin{equation}
\alpha=-\nu_1-2\nu_2+1,\;\; \beta=-\nu_1+\nu_2,\;\; \gamma=2\nu_1+\nu_2-1,
\end{equation}
for $k=0, 1;$ for $\phi(x)\in C_c^{\infty}(0, \infty)$ and $\tilde{\phi}(s):=\int\limits_0^\infty\phi(x)x^s\frac{dx}{x},$
set
\begin{equation}
\Phi_k(x):=\int\limits_{\Re s=\sigma}(\pi^3x)^{-s}\frac{
\Gamma\left(\frac{1+s+2k+\alpha}{2}\right)
\Gamma\left(\frac{1+s+2k+\beta}{2}\right)
\Gamma\left(\frac{1+s+2k+\gamma}{2}\right)}
{\Gamma\left(\frac{-s-\alpha}{2}\right)
\Gamma\left(\frac{-s-\beta}{2}\right)
\Gamma\left(\frac{-s-\gamma}{2}\right)}
\tilde{\phi}(-s-k)ds
\end{equation}
with $\sigma>\max\{-1-\Re\alpha, -1-\Re\beta, -1-\Re\gamma\},$
\begin{equation}
\Phi_{0, 1}^0(x)=\Phi_0(x)+\frac{\pi^{-3}c^3n}{m_1^2m_2i}\Phi_1(x)
\end{equation}
and
\begin{equation}
\Phi_{0, 1}^1(x)=\Phi_0(x)-\frac{\pi^{-3}c^3n}{m_1^2m_2i}\Phi_1(x),
\end{equation}
 we have the following Voronoi formula on ${\rm GL}(3):$
\begin{proposition} {\rm{([MS1], [GL])}} Let $\phi(x)\in C_c^{\infty}(0, \infty).\;$Let $A(n, m)$ denote the 
$(n, m)$-th Fourier coefficient of a Maass form for $SL(3, \mathbb{Z})$ as in 
(2.1).  Let
${a, \bar{a}, c} \in\mathbb{Z}$ with $c\neq 0, (a, c)=1,$
and $a\bar{a}\equiv 1 (\rm {mod} \;c).$ Then
%Then for any integer $\delta>0,$ $c\neq 0,$
 %$(a, c)=1$ and $a\bar{a}\equiv 1 (\rm{mod}\; c),$
we have
\begin{eqnarray}
\lefteqn{
\sum\limits_{m>0}A(n, m)e\left(\frac{m\bar{a}}{c}\right)\phi(m)\nonumber}\\
&&
=\frac{c\pi^{-\frac{5}{2}}}{4i}\sum\limits_{m_1|cn}\sum\limits
_{m_2>0}\frac{A(m_1, m_2)}{m_1m_2}S(na, m_2; n cm_1^{-1})\Phi_{0, 1}^0
\left(\frac{m_2m_1^2}{c^3n}\right)
\nonumber\\
&&\;\;\;\;+\frac{c\pi^{-\frac{5}{2}}}{4i}\sum\limits_{m_1|cn}\sum\limits
_{m_2>0}\frac{A(m_1, m_2)}{m_1m_2}S(n a, -m_2; n cm_1^{-1})
\Phi_{0, 1}^1
\left(\frac{m_2m_1^2}{c^3n}\right)
\nonumber,
\end{eqnarray}
where $S(a, b; c)$ is the Kloosterman sum defined as the above.
\end{proposition}
\section{$L$-functions}
\setcounter{equation}{0}  For each $u_j(z)$ of type $\frac{1}{2}+it_j$ in the orthonormal basis of even Maass forms for $SL(2,\mathbb{Z})$ 
with the normalized Fourier coefficients $a_j(n)$ as in (2.4), we associate the $L$-function $L(s, u_j)$ as in (1.1) which is entire and satisfies the following functional equation
\begin{equation}
\Lambda(s, u_j):=\pi^{-s}
\Gamma\left(\frac{s+it_j}{2}\right)\Gamma\left(\frac{s-it_j}{2}\right)L(s, u_j)=\Lambda(1-s, u_j).
\end{equation}
Using the functional equation (3.1) we shall represent the central values $L(\frac{1}{2}, u_j)$ by its partial sums of length about $O(|t_j|).$
To this need, we choose a function ([IK], pp. 98)
\begin{equation}
G(u)=\left(\cos\frac{\pi u}{A}\right)^{-A},
\end{equation}
we then have the following approximate functional equation (this has been worked out for general $L$-functions in [IK]):
\begin{lemma}
For any $u_j(z)$ of type $\frac{1}{2}+it_j$ in the orthonormal basis of even Maass forms for $SL(2, \mathbb{Z})$
\begin{equation}
L\left(\frac{1}{2}, u_j\right)=2\sum\limits_{l\geqslant 1}a_j(l)l^{-\frac{1}{2}}U(l, t_j)
\end{equation}
with
\begin{equation}
U(y, t)=\frac{1}{2\pi i}\int\limits_{\big(\frac{1}{2}\big)}y^{-u}G(u)
\frac{\gamma(\frac{1}{2}+u, t)}{\gamma(\frac{1}{2}, t)}\frac{du}{u}
\end{equation}
where $G(u)$ is defined by (3.2) and 
\begin{equation}
\gamma(u, t)=\pi^{-u}
\Gamma\left(\frac{u+it}{2}\right)
\Gamma\left(\frac{u-it}{2}\right).
\end{equation}
\end{lemma}
\noindent{\bf Proof.} See [IK] pp. 98. \mbox{$\Box$}\\
$U(y, t)$ has the following properties which effectively limit the terms in (3.3) with $l\ll |t_j|.$
\begin{lemma}
For $y, t>0,$ \\
1) (\rm {[IK]}, pp. 100) the derivatives of $U(y, t)$ with respect to $y$ satisfy 
$$
y^a\frac{\partial^a}{\partial y^a}U(y, t)\ll \left(1+\frac{y}{|t|}\right)^{-A},
$$
$$
y^a\frac{\partial^a}{\partial y^a}U(y, t)=\delta_a+O\left(\left(\frac{y}{|t|}\right)^\alpha\right),
$$
where
$0<\alpha\leqslant \frac{1}{6}, \delta_0=1, 0$ otherwise and the implied constants depend only on $\alpha, a$ and $A.$\\
2) if $1\leqslant y\ll t^{1+\varepsilon}, $ then we have the following asymptotic expansion as $t\rightarrow\infty$
\begin{equation}
U(y, t)=\frac{1}{2\pi i}\int\limits_{\big(\frac{1}{2}\big)}\left(\frac{t}{2\pi y}\right)^uG(u)\left[1+\frac{p_2(v)}{t}
+\frac{p_4(v)}{t^2}+ O\left(\frac{p_6(v)}{t^3}\right)\right]\frac{du}{u}+O(t^{-B}),
\end{equation}
where $v=\Im u,\; p_i(v)$ are polynomials of $v$ of degree $i$ and $B>0$ is arbitrarily large.
\end{lemma}
\noindent{\bf Proof.} 1) See [IK], pp. 100.\\
2) For $\Re u=\frac{1}{4}, \Im s\rightarrow \infty$ and $\Re s=\frac{1}{4},$ by Stirling's formula, we have
$$
\frac{\Gamma(s+u)}{\Gamma(s)}\ll |s|^{\frac{1}{4}}\exp\left(\frac{\pi}{2}|u|\right).
$$
It follows that for $\Re u=\frac{1}{2},$
$$
\frac{\Gamma\left(\frac{\frac{1}{2}+u+it}{2}\right)}{\Gamma\left(\frac{\frac{1}{2}+it}{2}\right)}\ll e^{\frac{\pi}{4}|u|}|t|^{\frac{1}{4}},
$$
$$
\frac{\Gamma\left(\frac{\frac{1}{2}+u-it}{2}\right)}{\Gamma\left(\frac{\frac{1}{2}-it}{2}\right)}\ll
 e^{\frac{\pi}{4}|u|}|t|^{\frac{1}{4}},
$$
hence
$$
\frac{\gamma(\frac{1}{2}+u, t)}{\gamma(\frac{1}{2}, t)}\ll e^{\frac{\pi}{2}|u|}|t|^{\frac{1}{2}},
$$
so
\begin{equation}
\frac{1}{2\pi i}\int\limits_{\substack{\left(\frac{1}{2}\right)\\|\Im u|\geqslant t^{\varepsilon}}}y^{-u}G(u)
\frac{\gamma(\frac{1}{2}+u, t)}{\gamma(\frac{1}{2}, t)}\frac{du}{u}
\ll |t|^{-B}
\end{equation}
for any large $B>0.$\\
By Stirling's formula
$$
\log \Gamma(s+c)=\left(s+c-\frac{1}{2}\right)\log s-s +\frac{1}{2}\log 2\pi +c_1s^{-1}+c_2s^{-2}+O\left(\frac{1}{|s|^3}\right)
$$
for any constant $c$ (the $c_{\nu}$'s are constants depending on $c$), as $|s|\rightarrow \infty$ uniformly for
 $|\arg s|\leqslant \pi-\varepsilon<\pi,$ one obtains that for $|u|\leqslant t^{\varepsilon},$
\begin{equation}
\frac
{\gamma(\frac{1}{2}+u, t)}{\gamma(\frac{1}{2}, t)}=\left(\frac{t}{2\pi}\right)^u\left[1+\frac{p_2(v)}{t}+\frac{p_4(v)}{t^2}
+O\left(\frac{p_6(v)}{t^2}\right)\right],
\end{equation}
where $p_i(v)$ are polynomials of $v$ of degree $i.$\\
Combining (3.4), (3.7) and (3.8) yield the conclusion of II). \mbox{$\Box$}\\
By Cauchy's inequality, Lemmas 3.1 and 3.2, we have 
\begin{equation}
L\left(\frac{1}{2}, u_j\right)\ll_{\varepsilon}|t_j|^{\frac{1}{2}+\varepsilon}
\end{equation}
where we used [Iw2] (pp. 130) and [HL], (3.9) is the convexity bound of $L(\frac{1}{2}, u_j).$
Correspondingly, to the Eisenstein series $E(z, s)$, we associate the \\
$L$-function
\begin{equation}
L(s, E)=\sum\limits_{n=1}^{\infty}\eta\left(n, \frac{1}{2}+ir\right)n^{-s}=\zeta(s-ir)\zeta(s+ir).
\end{equation}
It satisfies the functional equation (3.1) which can be verified directly using
the functional equation of $\zeta(s);$ so (3.3) becomes
\begin{equation}
\big|\zeta\left(\frac{1}{2}+ir\right)\big|^2=2\sum\limits_{l\geqslant 1}\left(\sum\limits_{ad=l}\left(\frac{a}{d}\right)^{ir}\right)l^{-\frac{1}{2}}U(l, r)
\end{equation}
where $U(l, r)$ is defined by (3.4).\\
Now let $f$ be a Maass form of type $(\nu_1, \nu_2)$ for $SL(3, \mathbb{Z}),$ the $L$-function $L(s, f)$ (see (1.2)) is entire and satisfies
the functional equation
$$
G_{\nu}(s)L(s, f)=\tilde{G}_{\nu}(1-s)L(1-s, \tilde{f})
$$
where
$$
G_{\nu}(s)=\pi^{\frac{-3s}{2}}
\Gamma\left(\frac{s+1-2\nu_1-\nu_2}{2}\right)
\Gamma\left(\frac{s+\nu_1-\nu_2}{2}\right)
\Gamma\left(\frac{s-1+\nu_1+2\nu_2}{2}\right),
$$
$$
G_{\tilde{\nu}}(s)=\pi^{\frac{-3s}{2}}
\Gamma\left(\frac{s+1-\nu_1-2\nu_2}{2}\right)
\Gamma\left(\frac{s-\nu_1+\nu_2}{2}\right)
\Gamma\left(\frac{s-1+2\nu_1+\nu_2}{2}\right),
$$
and $\tilde{f}$ is the dual Maass form.
The Rankin-Selberg $L$-function defined by
$$
L(s, f\times f):=\sum\limits_{m\geqslant 1}\sum\limits_{n\geqslant 1}\frac{|A(m, n)|^2}{(m^2n)^s}
$$
for $\Re s$ large has a meromorphic continuation to the whole plane with the only simple pole at $s=1.$ By a standard contour integration,
one shows that
\begin{equation}
\mathop{\sum\sum}_{m^2n\leqslant N}
|A(m, n)|^2\ll_f N.
\end{equation}
By Cauchy's inequality and (3.12), one derives that
\begin{equation}
\sum\limits_{n\leqslant N}|A(m, n)|\ll_f Nm.
\end{equation}
The Rankin-Selberg $L$-function of $f$ and $u_j$
$$
L(s, f\times u_j)=\sum\limits_{m\geqslant 1}\sum\limits_{n\geqslant 1}\frac{\bar{a}_j(n)A(m, n)}{(m^2n)^s}
$$
is entire and satisfies the functional equation
\begin{equation}
\Lambda(s, f\times u_j)=\Lambda(1-s, \tilde{f}\times u_j)
\end{equation}
where
\begin{eqnarray}
\lefteqn{
\Lambda(s, f\times u_j)=\pi^{-3s}
\Gamma\left(\frac{s-it_j-\alpha}{2}\right)
\Gamma\left(\frac{s-it_j-\beta}{2}\right)
\Gamma\left(\frac{s-it_j-\gamma}{2}\right)}\nonumber\\
&&\hspace{1.5cm}\cdotp
\Gamma\left(\frac{s+it_j-\alpha}{2}\right)
\Gamma\left(\frac{s+it_j-\beta}{2}\right)
\Gamma\left(\frac{s+it_j-\gamma}{2}\right)L(s, f\times u_j)
\nonumber
\end{eqnarray}
and 
\begin{eqnarray}
\lefteqn{
\Lambda(s, \tilde{f}\times u_j)=\pi^{-3s}
\Gamma\left(\frac{s+it_j+\alpha}{2}\right)
\Gamma\left(\frac{s+it_j+\beta}{2}\right)
\Gamma\left(\frac{s+it_j+\gamma}{2}\right)}\nonumber\\
&&\hspace{1.5cm}\cdotp
\Gamma\left(\frac{s-it_j+\alpha}{2}\right)
\Gamma\left(\frac{s-it_j+\beta}{2}\right)
\Gamma\left(\frac{s-it_j+\gamma}{2}\right)L(s, \tilde{f}\times u_j),
\nonumber
\end{eqnarray}
in the above,
$$
\alpha=-\nu_1-2\nu_2+1,\;\;\beta=-\nu_1+\nu_2,\;\;\gamma=2\nu_1+\nu_2-1.
$$
Set
\begin{equation}
F(u)=\left(\cos \frac{\pi u}{A}\right)^{-3A},
\end{equation}
\begin{equation}
V_1(y, t)=\frac{1}{2\pi i}\int\limits_{(3)}y^{-u}F(u)\frac{\gamma_1(\frac{1}{2}+u, t)}{\gamma_1(\frac{1}{2}, t)}
\frac{du}{u},
\end{equation}
\begin{equation}
V_2(y, t)=\frac{1}{2\pi i}\int\limits_{(3)}y^{-u}F(u)\frac{\gamma_2(\frac{1}{2}+u, t)}{\gamma_1(\frac{1}{2}, t)}
\frac{du}{u},
\end{equation}
note that one could move the line of integration in $V_1(y, t)$ and $V_2(y, t)$ to $\frac{1}{2}$ which is justified by Luo-Rudnick-Sarnak's 
bound on the Ramanujan conjecture $|\Re \alpha|, |\Re\beta|, |\Re \gamma|\leqslant \frac{1}{2}-\frac{1}{10}$ (see [LRS]), 
\begin{eqnarray}
\lefteqn{\;\;\;\;\;\;\;\;\;\;\;
\gamma_1(s, t)=\pi^{-3s}
\Gamma\left(\frac{s-it-\alpha}{2}\right)
\Gamma\left(\frac{s-it-\beta}{2}\right)
\Gamma\left(\frac{s-it-\gamma}{2}\right)}\\
&&\hspace{3.9cm}\cdotp
\Gamma\left(\frac{s+it-\alpha}{2}\right)
\Gamma\left(\frac{s+it-\beta}{2}\right)
\Gamma\left(\frac{s+it-\gamma}{2}\right),
\nonumber
\end{eqnarray}
\begin{eqnarray}
\lefteqn{\;\;\;\;\;\;\;\;\;\;\;\;
\gamma_2(s, t)=\pi^{-3s}
\Gamma\left(\frac{s-it+\alpha}{2}\right)
\Gamma\left(\frac{s-it+\beta}{2}\right)
\Gamma\left(\frac{s-it+\gamma}{2}\right)}\\
&&\hspace{3.9cm}\cdotp
\Gamma\left(\frac{s+it+\alpha}{2}\right)
\Gamma\left(\frac{s+it+\beta}{2}\right)
\Gamma\left(\frac{s+it+\gamma}{2}\right),
\nonumber
\end{eqnarray}
one has the following approximate functional equation for $L(s, f\times u_j):$
\begin{lemma}
For a Maass form $f$ of type $(\nu_1, \nu_2)$ for $SL(3, \mathbb{Z})$ and any $u_j(z)$ of type $\frac{1}{2}+it_j$ in the orthonormal
basis of even Hecke-Maass forms for $SL(2, \mathbb{Z}),$ 
we have

\begin{eqnarray}
\lefteqn{
L\left(\frac{1}{2}, f\times u_j\right)=\sum\limits_{m\geqslant 1}\sum\limits_{n\geqslant 1}
\frac{\bar{a}_j(n)A(m, n)}{(m^2n)^{\frac{1}{2}}}
V_1(m^2n, t_j)}\\
&&\hspace{2cm} +\sum\limits_{m\geqslant 1}\sum\limits_{n\geqslant 1}\frac{\bar{a}_j(n)A(n, m)}{(m^2n)^{\frac{1}{2}}}
V_2(m^2n, t_j).\nonumber
\end{eqnarray} 
\end{lemma}
\noindent{\bf Proof.} Following [IK] pp. 98, we consider the integral 
$$
I\left(\frac{1}{2}, f\times u_j\right)=\frac{1}{2\pi i}\int\limits_{(3)}\Lambda\left(\frac{1}{2}+u, f\times u_j\right) F(u)\frac{du}{u}.
$$
Moving the line of integration to $\Re u=-3$ and applying the functional equation, there yields 
$$
\Lambda \left(\frac{1}{2}, f\times u_j\right)=I\left(\frac{1}{2}, f\times u_j\right)+ 
I\left(\frac{1}{2}, \tilde{f}\times u_j\right)
$$
where $\Lambda(\frac{1}{2}, f\times u_j)$ comes from the simple pole of $u^{-1}F(u)$ at $u=0.$ By expanding into absolutely convergent 
Dirichlet series, we have
$$
I\left(\frac{1}{2}, f\times u_j\right)= \gamma_1\left(\frac{1}{2}, t_j\right)\sum\limits_{m\geqslant 1}\sum\limits_{n\geqslant 1}
\frac{\bar {a}_j(n)A(m, n)}{(m^2n)^{\frac{1}{2}}}
\frac{1}{2\pi i}\int\limits_{(3)}(m^2n)^{-u}F(u)\frac{\gamma_1(\frac{1}{2}+u, t_j)}{\gamma_1(\frac{1}{2}, t_j)}
\frac{du}{u}.$$ Similarly, 
$$
I\left(\frac{1}{2}, \tilde{f}\times u_j\right)= \gamma_1\left(\frac{1}{2}, t_j\right)\sum\limits_{m\geqslant 1}\sum\limits_{n\geqslant 1}
\frac{\bar {a}_j(n)A(n, m)}{(m^2n)^{\frac{1}{2}}}
\frac{1}{2\pi i}\int\limits_{(3)}(m^2n)^{-u}F(u)\frac{\gamma_2(\frac{1}{2}+u, t_j)}{\gamma_1(\frac{1}{2}, t_j)}
\frac{du}{u}.$$ 
Combining them and dividing both sides by $\gamma_1(\frac{1}{2}, t_j),$ one finishes the proof of the lemma.
\mbox{$\Box$}\\
$V_1(y, t)$ and $V_2(y, t)$ have the following properties which effectively limit the terms in (3.20) with $m^2n\ll |t_j|^3.$
\begin{lemma}For $y, t>0, i=1, 2,$\\
1) the derivatives of $V_i(y, t)$ with respect to $y$ satisfy 
$$
y^a\frac{\partial^a}{\partial y^a}V_i(y, t)\ll \left(1+\frac{y}{|t|^3}\right)^{-A},
$$
$$
y^a\frac{\partial^a}{\partial y^a}V_i(y, t)=\delta_a+O\left(\left(\frac{y}{|t|^3}\right)^{c}\right)
$$
%\end{document}
where
$0<c\leqslant \frac{1}{3}\text{min} \{
\frac{1}{2}-\Re \alpha, \frac{1}{2}-\Re \beta, \frac{1}{2}-\Re \gamma\},\;
\delta_0=1, 0$ otherwise and the implied constants depend only on $\alpha, a, A$ and $f.$\\
2) if $1\leqslant y\ll t^{3+\varepsilon},$ then as $t\rightarrow\infty,$ we have
\begin{equation}
V_i(y, t)=\frac{1}{2\pi i}
\int\limits_{\left(\frac{1}{2}\right)}
\left(\frac{t^3}{8\pi^3 y}\right)^uF(u)\left[1+\frac{p_2(v)}{t}
+\frac{p_4(v)}{t^2}+O\left(\frac{p_6(v)}{t^3}\right)\right]\frac{du}{u}+O\left(t^{-B}\right)
\end{equation}
where $v=\Im u,$ $p_i(v)$ are polynomials of $v$ of degree $i$ and $B$ is arbitrarily large.
\end{lemma}
\noindent{\bf Proof.} 1) See [IK], pp. 100.\\
2) Similar to the proof of Lemma 3.2 2). \mbox{$\Box$}\\
%\end{document}
By Lemmas 3.3 and 3.4,
\begin{eqnarray}
\lefteqn{
L\left(\frac{1}{2}, f\times u_j\right)\ll 
\mathop{\sum\sum}_{m^2n\leqslant |t_j|^{3+\varepsilon}}
\frac{|a_j(n)A(m, n)|}{(m^2n)^{\frac{1}{2}}}
|t_j|^{\varepsilon}}\nonumber\\
&&\hspace{3cm}+\mathop{\sum\sum}\limits_{m^2n\leqslant |t_j|^{3+\varepsilon}}
\frac{|a_j(n)A(n, m)|}{(m^2n)^{\frac{1}{2}}}
|t_j|^{\varepsilon}\nonumber.
\end{eqnarray}
Furthurmore, applying Cauchy's inequality, [Iw2] (pp. 130), [HL] and (3.12), we have
\begin{equation}
\mathop{\sum\sum}_{m^2n\leqslant |t_j|^{3+\varepsilon}}
\frac{|a_j(n)A(m, n)|}{(m^2n)^{\frac{1}{2}}}\ll_{f, \varepsilon} |t_j|^{\frac{3}{2}+\varepsilon}.
\end{equation}
%\end{document}
Similarly, one can prove that
\begin{equation}
\mathop{\sum\sum}_{m^2n\leqslant |t_j|^{3+\varepsilon}}
\frac{|a_j(n)A(n, m)|}{(m^2n)^{\frac{1}{2}}}\ll_{f, \varepsilon} |t_j|^{\frac{3}{2}+\varepsilon}.
\end{equation}
%\end{document}
Combining (3.22) and (3.23), we obtain the convexity bound
\begin{equation}
L\left(\frac{1}{2}, f\times u_j\right)\ll_{f, \varepsilon}|t_j|^{\frac{3}{2}+\varepsilon}.
\end{equation}
%\end{document}
It follows from the convexity bound (3.9) of $L\left(\frac{1}{2}, u_j\right)$ and Weyl's law, the contribution to (1.4) from the error
term of (3.21) is bounded by
$$
{\sum_j}^{'} e^{-\frac{t_j^2}{T^2}}t_j^{\frac{1}{2}+\varepsilon}\mathop{\sum\sum}_{m^2n\leqslant |t_j|^{3+\varepsilon}}
\frac{|a_j(n)A(m, n)|}{m^2n}\frac{1}{t_j^{\frac{3}{2}}}
=O(T^{1+\varepsilon})
$$
where we also used [Iw2] (pp. 130), [HL] and (3.12).\\
%\end{document}
Similarly, the contribution to (1.4) from the error term of (3.6) is bounded by $O(T^{1+\varepsilon}). $\\
From now on, we only consider the leading terms in (3.6) and (3.21) since all the other terms can be treated similarly.\\
To the Maass form $f$ of type $(\nu_1, \nu_2)$ for $SL(3, \mathbb{Z})$ and the Eisenstein series $E\left(z, \frac{1}{2}+ir\right)$
we associate the $L$-function
$$
L(s, f\times E):=\sum\limits_{m\geqslant 1}\sum\limits_{n\geqslant 1}\frac{\bar{\eta}(n, \frac{1}{2}+ir)A(m, n)}{(m^2n)^s}
.$$
%\end{document}
By looking at the Euler products
$$
L(s, f)=\sum\limits_{m\geqslant 1}\frac{A(1, m)}{m^s}=
\prod\limits_p\prod\limits_{i=1}^3(1-\beta_{p, i}p^{-s})^{-1},
$$
$$
L(s, E)=\sum\limits_{n\geqslant 1}\eta\left(n, \frac{1}{2}+ir\right)n^{-s}=\prod\limits_p(1-p^{-s+ir})^{-1}(1-p^{-s-ir})^{-1},
$$
%\end{document}
one derives that
\begin{eqnarray}
\lefteqn{
L(s, f\times E)=\prod\limits_p\prod\limits_{k=1}^3
(1-\beta_{p, k}p^{ir-s})^{-1}(1-\beta_{p, k}p^{-ir-s})^{-1}}\nonumber\\
&&\hspace{5cm}=L(s-ir, f)L(s+ir, f).\nonumber
\end{eqnarray}
%\end{document}
It yields that 
\begin{equation}
L\left(\frac{1}{2}, f\times E\right)=L\left(\frac{1}{2}-ir, f\right)L\left(\frac{1}{2}+ir, f\right).
\end{equation}
%\end{document}
This satisfies the functional equation (3.14) which can also be verified directly using the functional equation of $L(s, f).$
So (3.20) becomes
\begin{eqnarray}
\lefteqn{\;\;\;\;\;\;\;\;\;
L\left(\frac{1}{2}, f\times E\right)
=
\sum\limits_{m\geqslant 1}\sum\limits_{n\geqslant 1}
\frac{\left(\sum\limits_{ad=n}\left(\frac{a}{d}\right)^{-ir}\right)A(m, n)}
{(m^2n)^{\frac{1}{2}}}
V_1(m^2n, r)}
\\
&&\hspace{4.2cm} +\sum\limits_{m\geqslant 1}\sum\limits_{n\geqslant 1}
\frac{\left(\sum\limits_{ad=n}\left(\frac{a}{d}\right)^{-ir}\right)A(n, m)}{(m^2n)^{\frac{1}{2}}}
V_2(m^2n, r).
\nonumber
\end{eqnarray}
%\end{document}
In using the Kuznetsov formula, we need also consider the continuous spectrum $E(z, s).$ We are led to prove the following proposition 
in order to prove our main theorem - Theorem 1.1:
%\newpage
\begin{proposition}
Let $f$ be a fixed Hecke-Maass form for $SL(3, \mathbb{Z}),$ $u_j$ an \\orthonormal basis of even Hecke-Maass forms for 
$SL(2, \mathbb{Z}),$
we have
\begin{eqnarray}
\lefteqn{{\sum_j}^{'} e^{-\frac{t_j^2}{T^2}}L\left(\frac{1}{2}, f\times u_j\right)L\left(\frac{1}{2}, u_j\right)}\\
&&+\int\limits_{-\infty}^{\infty}e^{-\frac{t^2}{T^2}}
\frac{L\left(\frac{1}{2}+it, f\right)L\left(\frac{1}{2}-it, f\right)|\zeta\left(\frac{1}{2}
+it\right)\big|^2}{\big|\zeta(1+2it)\big|^2}dt\nonumber\\
&&\hspace{2cm}=\frac{12T^2}{\pi^3}L(1, f)L(1, \tilde{f})+O_{\varepsilon, f}\left(T^{\frac{11}{6}+\varepsilon}\right).\nonumber
\end{eqnarray}
\end{proposition}
\noindent{\bf Remarks.} 1. Because of (3.10) and (3.25), one can see that the integral in the above is the 
continuous analogue of the discrete part. Actually, the contribution from the integral on the left of (3.27) is small. 
Indeed, by the well known bounds [Ti]
$$
\zeta(1+2it)\gg \log (1+2|t|)^{-1},
$$
$$
\zeta\left(\frac{1}{2}+it\right)\ll (|t|+1)^{\frac{1}{6}+\varepsilon}
$$
and 
$$
\int\limits_{-\infty}^{\infty}
e^{-\frac{t^2}{T^2}}\left|L\left(\frac{1}{2}+it, f\right)\right|^2dt
\ll T^{\frac{3}{2}+\varepsilon}
$$
which is a direct consequence of the approximate functional equation of $L(s, f)$ (see [IK], pp. 98-100), one derives that
\begin{equation}
\int\limits_{-\infty}^{\infty}e^{-\frac{t^2}{T^2}}\frac{L\left(\frac{1}{2}+it, f\right)L\left(\frac{1}{2}-it, f\right)|\zeta\left(\frac{1}{2}
+it\right)\big|^2}{\big|\zeta(1+2it)|^2}dt\ll_{\varepsilon, f}T^{\frac{11}{6}+\varepsilon}
\end{equation}
which is admissible with the error term in Theorem 1.1.\\
Let $\Omega(x)$ be a smooth function compactly supported on $\left[\frac{1}{2}, 2T^{1+\varepsilon}\right]$ 
with $\Omega=1$ on
$[1, T^{1+\varepsilon}]$ and $0\leqslant \Omega\leqslant 1$ on $\left[\frac{1}{2}, 2T^{1+\varepsilon}\right];$ 
$k(x)$ be a smooth
function compactly supported on $\left[\frac{1}{2}, 2T^{3+\varepsilon}\right]$ with $k(x)=1$ on
$[1, T^{3+\varepsilon}]$ and $0\leqslant k(x)\leqslant 1$ on $\left[\frac{1}{2}, 2T^{3+\varepsilon}\right],$ then (3.3) and (3.20)
yield that 
\begin{eqnarray}
\lefteqn{\hspace{0.7cm}
{\sum\limits_j}^{'}
 e^{-\frac{t_j^2}{T^2}}L\left(\frac{1}{2}, f\times u_j\right)L\left(\frac{1}{2}, u_j\right)}\nonumber\\
&&=2{\sum\limits_j}^{'}e^{-\frac{t_j^2}{T^2}}\sum\limits_{l\geqslant 1}\sum\limits_{m\geqslant 1}\sum\limits_{n\geqslant 1}
\frac{a_j(l)\bar{a}_j(n)A(m, n)}
{(m^2nl)^{\frac{1}{2}}}U(l, t_j)V_1(m^2n, t_j)\Omega(l)k(m^2n)\nonumber\\
&&\label{(3.29)}\\
&&+2{\sum\limits_j}^{'}e^{-\frac{t_j^2}{T^2}}\sum\limits_{l\geqslant 1}\sum\limits_{m\geqslant 1}\sum\limits_{n\geqslant 1}
\frac{a_j(l)\bar{a}_j(n)A(n, m)}
{(m^2nl)^{\frac{1}{2}}}U(l, t_j)V_2(m^2n, t_j)\Omega(l)k(m^2n)\nonumber\\
&&\hspace{9.5cm}+O_{f, \varepsilon}(T^{-B})\nonumber
\end{eqnarray}
where $B>0$ is arbitrarily large.\\
Next we transform the main term in (3.29) by the Kuznetsov formula (2.5) into $\Delta+\frac{1}{2}N\Delta,$
where
\begin{eqnarray}
\lefteqn{\Delta=
\sum\limits_{l\geqslant 1}\sum\limits_{m\geqslant 1}\sum\limits_{n\geqslant 1}\frac{A(m, n)}
{(m^2nl)^{\frac{1}{2}}}\Omega(l)k(m^2n)\delta(n, l)H_1}\\
&&\hspace{2cm}+\sum\limits_{l\geqslant 1}\sum\limits_{m\geqslant 1}\sum\limits_{n\geqslant 1}\frac{A(n, m)}
{(m^2nl)^{\frac{1}{2}}}\Omega(l)k(m^2n)\delta(n, l)H_2\nonumber
\end{eqnarray}
is the diagonal term, \\
$N\Delta=N\Delta^{1, 1}+N\Delta^{1, 2}+N\Delta^{2, 1}+N\Delta^{2, 2}$ is the nondiagonal term with
\begin{eqnarray}
\lefteqn{N\Delta^{1,1}=2\sum\limits_{l\geqslant 1}\sum\limits_{m\geqslant 1}\sum\limits_{n\geqslant 1}\frac{A(m, n)}{(m^2nl)^{\frac{1}{2}}}
\Omega(l)k(m^2n)}\\
&&\hspace{3cm}\cdotp\sum\limits_{c>0}c^{-1}S(n, l; c)H_1^+\left(\frac{2\sqrt{nl}}{c}\right),\nonumber
\end{eqnarray}
\begin{eqnarray}
\lefteqn{N\Delta^{1,2}=2\sum\limits_{l\geqslant 1}\sum\limits_{m\geqslant 1}\sum\limits_{n\geqslant 1}\frac{A(n, m)}{(m^2nl)^{\frac{1}{2}}}
\Omega(l)k(m^2n)}\\
&&\hspace{3cm}\cdotp\sum\limits_{c>0}c^{-1}S(n, l; c)H_1^+\left(\frac{2\sqrt{nl}}{c}\right),\nonumber
\end{eqnarray}
\begin{eqnarray}
\lefteqn{N\Delta^{2,1}=2\sum\limits_{l\geqslant 1}\sum\limits_{m\geqslant 1}\sum\limits_{n\geqslant 1}
\frac{A(m, n)}{(m^2nl)^{\frac{1}{2}}}
\Omega(l)k(m^2n)}\\
&&\hspace{3cm}\cdotp\sum\limits_{c>0}c^{-1}S(-n, l; c)H_1^-\left(\frac{2\sqrt{nl}}{c}\right),\nonumber
\end{eqnarray}
\begin{eqnarray}
\lefteqn{N\Delta^{2,2}=2\sum\limits_{l\geqslant 1}\sum\limits_{m\geqslant 1}\sum\limits_{n\geqslant 1}\frac{A(n, m)}{(m^2nl)^{\frac{1}{2}}}
\Omega(l)k(m^2n)}\\
&&\hspace{3cm}\cdotp\sum\limits_{c>0}c^{-1}S(-n, l; c)H_1^-\left(\frac{2\sqrt{nl}}{c}\right),\nonumber
\end{eqnarray}
\begin{equation}
H_1=\frac{1}{\pi}\int\limits_{-\infty}^{\infty}e^{-\frac{t^2}{T^2}}U(n, t)V_1(m^2n, t)\tanh (\pi t) tdt,
\end{equation}
\begin{equation}
H_1^+(x)=2i\int\limits_{-\infty}^{\infty}J_{2it}(2\pi x)\frac{e^{-\frac{t^2}{T^2}}U(l, t)V_1(m^2n, t)t}{\cosh \pi t}dt,
\end{equation}
\begin{equation}
H_1^-(x)=\frac{4}{\pi}\int\limits_{-\infty}^\infty K_{2it}(2\pi x)\sinh \pi te^{-\frac{t^2}{T^2}}U(l, t)V_1(m^2n, t)tdt,
\end{equation}
\begin{equation}
H_2=\frac{1}{\pi}\int\limits_{-\infty}^{\infty}e^{-\frac{t^2}{T^2}}U(n, t)V_2(m^2n, t)\tanh (\pi t) tdt,
\end{equation}
\begin{equation}
H_2^+(x)=2i\int\limits_{-\infty}^\infty J_{2it}(2\pi x)\frac{e^{-\frac{t^2}{T^2}}U(l, t)V_2(m^2n, t)t}{\cosh \pi t}dt,
\end{equation}
\begin{equation}
H_2^-(x)=\frac{4}{\pi}\int\limits_{-\infty}^{\infty}
K_{2it}(2\pi x)\sinh \pi t e^{-\frac{t^2}{T^2}}U(l, t)V_2(m^2n, t)tdt.
\end{equation}
\section{Evaluation of the diagonal terms}
\setcounter{equation}{0} 
In this section, we will estimate the contribution from the diagonal term $\Delta$ which is defined by (3.30).\\
Write $\Delta$ as $\Delta_1+\Delta_2,$ where
\begin{equation}
\Delta_1=\sum\limits_{n\geqslant 1}\sum\limits_{m\geqslant 1}\frac{A(m, n)}{mn}\Omega(n)k(m^2n)H_1
\end{equation}
and
\begin{equation}
\Delta_2=\sum\limits_{n\geqslant 1}\sum\limits_{m\geqslant 1}\frac{A(n, m)}{mn}\Omega(n)k(m^2n)H_2.
\end{equation}
Clearly
\begin{equation}
\Delta_1=\sum\limits_{n\geqslant 1}\sum\limits_{m\geqslant 1}\frac{A(m, n)}{mn}H_1+O(T^{-B})
\end{equation}
where $B>0$ is arbitrarily large.\\
Let's first consider 
\begin{equation}
\Delta_1^*:=\sum\limits_{n\geqslant 1}\sum\limits_{m\geqslant 1}\frac{A(m, n)}{mn}U(n, t)V_1(m^2n, t).
\end{equation}
Set
$$
\tilde{U}(s, t):=\int\limits_0^\infty U(x, t)x^s\frac{dx}{x}
$$
which is equal to $\frac{G(s)}{s}\frac{\gamma\left(\frac{1}{2}+s, t\right)}{\gamma\left(\frac{1}{2}, t\right)}
$ by (3.4),
set
$$
\tilde{V_1}(s_1, t):=\int\limits_0^\infty V_1(x, t)x^{s_1}\frac{dx}{x}
$$
which is equal to $\frac{F(s_1)}{s_1}\frac{\gamma_1\left(\frac{1}{2}+s_1, t\right)}{\gamma_1\left(\frac{1}{2}, t\right)}
$ by (3.16).
The Mellin inversion formula yields that
$$
V_1(y, t)=
\frac{1}{2\pi i}\int\limits_{(\sigma_1)}\tilde{V_1}(s_1, t)y^{-s_1}ds_1
$$
and 
$$
U(x, t)=\frac{1}{2\pi i}\int\limits_{(\sigma)}\tilde{U}(s, t)x^{-s}ds
$$
with $\sigma>-\frac{1}{2}$ and $\sigma_1>-\frac{1}{10}$ which is justified by Luo-Rudnick-Sarnak's bound on the generalized Ramanujan conjecture [LRS].\\
Due to Bump [Bu], we know that
$$
\sum\limits_{m\geqslant 1}\sum\limits_{n\geqslant 1}\frac{A(m, n)}{m^{s+1}n^{w+1}}=
\frac{L(s+1, \tilde{f})L(w+1, f)}{\zeta(s+w+2)}
,$$
therefore
\begin{eqnarray}
\lefteqn{
\Delta_1^*=
\frac{1}{(2\pi i)^2}
\int\limits_{(3)}\int\limits_{(3)}
s_1\tilde{V_1}(s_1, t)s\tilde{U}(s, t)
\sum\limits_{m\geqslant 1}\sum\limits_{n\geqslant 1}\frac{A(m, n)}{mn}n^{-s}
(m^2n)^{-s_1}\frac{ds}{s}\frac{ds_1}{s_1}}\nonumber\\
&&=\frac{1}{(2\pi i)^2}
\int\limits_{(3)}\int\limits_{(3)}
s_1\tilde{V_1}(s_1, t)s\tilde{U}(s, t)
\frac{L(2s_1+1, \tilde{f})L(s_1+s+1, f)}{\zeta(3s_1+s+2)}\frac{ds}{s}\frac{ds_1}{s_1}.
\nonumber
\end{eqnarray}
Moving the line of integration to $\sigma=-\frac{1}{4}$ and $\sigma_1=-\frac{1}{11},$ picking up a pole at $(0, 0),$ by the Residue theorem,
we have
\begin{eqnarray}
\lefteqn{
\Delta_1^*=\frac{L(1, \tilde{f})L(1, f)}{\zeta(2)}+\frac{1}{2\pi i}\int\limits_{\left(-\frac{1}{4}\right)}\frac{L(1, \tilde{f})L(s+1, f)}{\zeta(s+2)}
s\tilde{U}(s, t)\frac{ds}{s}}\nonumber\\
&&+\left(\frac{1}{2\pi i}\right)^2\int\limits_{\left(-\frac{1}{4}\right)}\int\limits_{\left(-\frac{1}{11}\right)}
\frac{L(2s_1+1, \tilde{f})L(s_1+s+1, f)}{\zeta(3s_1+s+2)}
\tilde{V_1}(s_1, t)\tilde{U}(s, t)dsds_1\nonumber\\
&&\hspace{4cm}+\frac{1}{2\pi i}\int\limits_{\left(-\frac{1}{11}\right)}\frac{L(2s_1+1, \tilde{f})L(s_1+1, f)}{\zeta(3s_1+2)}
\tilde{V_1}(s_1, t)ds_1\nonumber\\
&&\hspace{6cm}=\frac{L(1, \tilde{f})L(1, f)}{\zeta(2)}+O\left((|t|+1)^{-\frac{1}{4}}\right).
\nonumber
\end{eqnarray}
Thus
\begin{equation}
\Delta_1=\frac{6T^2}{\pi^3}L(1, f)L(1, \tilde{f})+O_{f, \varepsilon}\left(T^{\frac{7}{4}+\varepsilon}\right).
\end{equation}
Similarly
\begin{equation}
\Delta_2=\frac{6T^2}{\pi^3}L(1, f)L(1, \tilde{f})+O_{f, \varepsilon}\left(T^{\frac{7}{4}+\varepsilon}\right).
\end{equation}
Combining (4.5) and (4.6), we obtain
\begin{equation}
\Delta=\frac{12T^2}{\pi^3}L(1, f)L(1, \tilde{f})+O_{f, \varepsilon}\left(T^{\frac{7}{4}+\varepsilon}\right).
\end{equation}
\section{Sums of Kloosterman sums - large $c$}
\setcounter{equation}{0} 
From now on, we will start to show the contribution from sums of Kloosterman sums - the nondiagonal terms is small
(recall (3.31)-(3.34)). For simplicity, we only treat $N\Delta^{1,1}$ later on since $N\Delta^{1, 2}, N\Delta^{2, 1}$
and $N\Delta^{2, 2}$ can be estimated similarly.\\
Let $g, h$ be smooth functions supported on $[1, 2],$ we apply smooth partitions of unity to $l$ and $m^2n,$
$$
1=\sum\limits_{u=-\infty}^{\infty}g\left(\frac{x}{2^u}\right),\;\;
 1=\sum\limits_{v=-\infty}^{\infty}h\left(\frac{x}{2^v}\right)
$$
then
\begin{eqnarray}
\lefteqn{
N\Delta^{1, 1}=2\sum\limits_{N_1}\sum\limits_{N_2}
\sum\limits_{l\geqslant 1} 
\sum\limits_{m\geqslant 1} 
\sum\limits_{n\geqslant 1}
\frac{A(m, n)}{(m^2nl)^{\frac{1}{2}}}\Omega(l)k(m^2n)}\nonumber\\
&&\cdotp g\left(\frac{l}{N_2}\right)h\left(\frac{m^2n}{N_1}\right)\sum\limits_{c>0}c^{-1}S(n, l; c)
H_1^+\left(\frac{2\sqrt{nl}}{c}\right)
\nonumber
\end{eqnarray}
where $N_2=2^u$ and $N_1=2^v.$\\
Since $\Omega(x)$ limits the $l$-sum to $l\ll T^{1+\varepsilon}$ and $k(x)$ limits $m^2n$ to \\
$m^2n\ll T^{3+\varepsilon}, N_2\ll T^{1+\varepsilon}$
and $N_1\ll T^{3+\varepsilon}.$ For fixed $m,$ we split the $c$-sum into three ranges:\\
I) $c\leqslant T^{\frac{11}{9}+\varepsilon}m^{-1};$ \;\;\\
II) $T^{\frac{11}{9}+\varepsilon}m^{-1}\leqslant c\leqslant Cm^{-1};$\;\;\\
III) $c\geqslant Cm^{-1}$ with $C=T^{\frac{11}{9}+\varepsilon}+\sqrt{N_1N_2}.$\\
In this section, we will study the third case.\\
Let
\begin{equation}
U(l, -\sigma i+y)=\frac{1}{2\pi i}\int\limits_{(A)}l^{-u}G(u)\frac{\gamma\left(\frac{1}{2}+u, -\sigma i+y\right)}
{\gamma\left(\frac{1}{2}, -\sigma i+y\right)}\frac{du}{u}
\end{equation}
and
\begin{equation}
V_j(m^2n, -\sigma i+y)=\frac{1}{2\pi i}\int\limits_{(B)}{(m^2n)}^{-u}F(u)\frac{\gamma_j\left(\frac{1}{2}+u, -\sigma i+y\right)}
{\gamma_1\left(\frac{1}{2}, -\sigma i+y\right)}\frac{du}{u}
\end{equation}
for $j=1, 2, A>\sigma-\frac{1}{2}$ and $B>\max \{
\sigma+|\Re \alpha|-\frac{1}{2}, \sigma+|\Re \beta|-\frac{1}{2}, \sigma+|\Re \gamma|-\frac{1}{2}\},$ where $\gamma(u, s)$ is 
defined by (3.5), $\gamma_1(u, s)$ and $\gamma_2(u, s)$ are defined by (3.18) and (3.19) respectively. Recall $G(u)$ is defined
by (3.2) and $F(u)$ is defined by (3.15).
\\
By Stirling's formula, one derives that
\begin{equation}
U(l, -\sigma i+y)\ll_{\sigma}\left(\frac{|y|}{l}\right)^A,
\end{equation}
\begin{equation}
V_j(m^2n, -\sigma i+y)\ll_{\sigma}\left(\frac{|y|^3}{m^2n}\right)^B.
\end{equation}
Recall $H_1^+(x)$ is defined by (3.36). Moving the line of integration to $\Im t=-\sigma,$ then $H_1^+(x)$ becomes
\begin{equation}
2i\!\!\int\limits_{-\infty}^{\infty}\!J_{2iy+2\sigma}(2\pi x)\frac{e^{-\frac{(-\sigma i+y)^2}{T^2}}U(l, -\sigma i+y)V_1(m^2n, -\sigma i+y)(-\sigma i+y)}{\cosh \pi (-\sigma i+y)}
dy.
\end{equation}
For $0<x<1,$ using the bound 
$$
J_{2iy+2\sigma}(2\pi x)\ll 
\left\{
\begin{array}{lll}
x^{2\sigma}e^{\pi |y|}|y|^{-2\sigma}& \mbox{} &\text{if} \;\;\;|y|\geqslant 1 \\
x^{2\sigma}&\mbox{}&  \text{if}\;\;\; |y|\leqslant 1 \end{array} \right.
$$
and (5.5), we have (taking $\sigma=A=B$)
$$
H_1^+(x)\ll x^{2\sigma}T^{2\sigma+2}l^{-\sigma}(m^2n)^{-\sigma}.
$$
By the above bound, (3.13) and Weil's bound for Kloosterman sums
\begin{equation}
S(n, l; c)\ll c^{\frac{1}{2}}(n, l, c)^{\frac{1}{2}}\tau(c),
\end{equation} we have 
\begin{eqnarray}
\lefteqn{
\sum\limits_{l\geqslant 1}\sum\limits_{m\geqslant 1}\sum\limits_{n\geqslant 1}
\frac{A(m, n)}{(m^2nl)^{\frac{1}{2}}}
\Omega(l)k(m^2n)g\left(\frac{l}{N_2}\right)
h\left(\frac{m^2n}{N_1}\right)}\nonumber\\
&&\hspace{2cm}\cdotp\sum\limits_{c\geqslant Cm^{-1}}c^{-1}S(n, l; c)
H_1^+\left(\frac{2\sqrt{nl}}{c}\right)
\nonumber\\
&&\hspace{5.5cm}\ll T^{-\frac{4}{9}\sigma+\frac{110}{18}}\nonumber
\end{eqnarray}
which is negligible since $\sigma$ can be very large. 
\section{Sums of Kloosterman sums - small $c$: Part I}
\setcounter{equation}{0} 
In the following two sections, we will estimate the contribution from \\$c\leqslant T^{\frac{11}{9}+\varepsilon}m^{-1}$ 
in $N_\Delta^{1, 1}$. The Voronoi formula
on ${\rm GL}(3)$ will be used.\\
The $J$-Bessel function has an integral representation ([GR], pp. 902):
$$
J_{2it}(2\pi x)=\frac{2}{\pi}\int\limits_0^\infty\sin (2\pi x\cosh \zeta-i\pi t)\cos 2t\zeta d\zeta
$$
from which one derives
$$
\frac{J_{2it}(2\pi x)-J_{-2it}(2\pi x)}{\cosh \pi t}=-\frac{2i}{\pi}\frac{\sinh \pi t}{\cosh \pi t}\int\limits_{-\infty}^{\infty}
\cos (2\pi x\cosh \zeta)\cos 2t\zeta d\zeta.
$$
For $x=\frac{2\sqrt{nl}}{c}\geqslant T^{-\frac{11}{9}-\varepsilon}$ and $|t|\leqslant T^{1+\varepsilon},$ by partial 
integration once, we have
\begin{eqnarray}
\lefteqn{
\frac{J_{2it}(2\pi x)-J_{-2it}(2\pi x)}{\cosh \pi t}}\\
%=\frac{2i}{\pi}\frac{\sinh \pi t}{\cosh \pi t}\int\limits_{-T^{\varepsilon}}^{T^{\varepsilon}}
%\cos(2\pi c\cosh\zeta)\cos 2t\zeta d\zeta}\\
&&\hspace{0cm}=\frac{2i}{\pi}\frac{\sinh \pi t}{\cosh \pi t}\int\limits_{-T^{\varepsilon}}^{T^{\varepsilon}}
\cos(2\pi x\cosh\zeta)\cos 2t\zeta d\zeta+O(T^{-B})\nonumber
\end{eqnarray}
for any $B>0.$\\
Combining (6.1) with the definition of $H_1^+(x)$ (recall (3.36)), one obtains 
\begin{eqnarray}
\lefteqn{
H_1^+(x)=-\frac{2}{\pi}\int\limits_{-T^{1+\varepsilon}}^{T^{1+\varepsilon}}\int\limits_{-T^{\varepsilon}}^{T^\varepsilon}
\frac{\sinh \pi t}{\cosh \pi t}e^{-\frac{t^2}{T^2}}U(l, t)V_1(m^2n, t)t}\\
&&\hspace{3cm}\cdotp \cos (2\pi x\cosh \zeta)\cos(2t\zeta)d\zeta dt+O(T^{-B}).
\nonumber
\end{eqnarray}
For convinence, we apply a smooth partition of unity to the variable $t$
$$
1=\sum\limits_{\alpha=-\infty}^\infty\eta\left(\frac{t}{2^{\alpha}}\right)
$$
where $\eta(x)$ is a smooth function compactly supported on $[1, 2]$, then
\begin{eqnarray}
\lefteqn{\;\;
H_1^+(x)=-\frac{2}{\pi}\sum\limits_{T_0}\int\limits_{-T^{1+\varepsilon}}^{T^{1+\varepsilon}}\int\limits_{-T^{\varepsilon}}^{T^\varepsilon}
\frac{\sinh \pi t}{\cosh \pi t}e^{-\frac{t^2}{T^2}}U(l, t)V_1(m^2n, t)t}\\
&&\hspace{2cm}\cdotp \cos (2\pi x\cosh \zeta)\cos(2t\zeta)\eta\left(\frac{t}{T_0}\right)d\zeta dt+O(T^{-B})
\nonumber
\end{eqnarray}
where $T_0=2^{\alpha}\ll T^{1+\varepsilon}.$\\
There are two cases:\\
I) If $N_2\geqslant T_0^{1+\varepsilon}$ or $N_1\geqslant T_0^{3+\varepsilon},$ then due to the rapid decay of $U(l, t)$ and $V_1(m^2n, t)$ (see Lemma 3.2 and Lemma 3.4), the contribution from such
terms is negligible.\\
II) If $N_2\leqslant T_0^{1+\varepsilon}$ and $N_1\leqslant T_0^{3+\varepsilon},$ we apply the asymptotic expansion of $U(l, t)$ and $V_1(m^2n, t)$
(see Lemma 3.2 and Lemma 3.4). If $|\zeta|\geqslant T_0^{-1+\varepsilon},$ applying partial integrations to the $t$-integral many times,
one shows that its contribution is negligible. So next we only consider the case $|\zeta|\leqslant T_0^{-1+\varepsilon}.$ We cut the 
$\zeta$-integral smoothly by inserting a smooth factor $\omega\left(\frac{\zeta}{T_0^{-1+\varepsilon}}\right),$ where $\omega(x)$ is compactly
supported on $[-2, 2]$ and equals 1 on $[-1, 1].$ 
Set
$$
W(t):=\int\limits_{-\infty}^{\infty}\omega\left(\frac{\zeta}{T_0^{-1+\varepsilon}}\right)\cos (2\pi x\cosh \zeta)\cos (2t\zeta)d\zeta
$$
then by partial integration twice, we have
$$
W(t)=W_1(t)+W_2(t)+O(T_0^{-B})
$$
with
$$
W_1(t)=\frac{\pi x}{2t^2}
\int\limits_{-\infty}^\infty\omega\left(\frac{\zeta}{T_0^{-1+\varepsilon}}\right)\cosh \zeta\sin(2\pi x\cosh \zeta)\cos 2t\zeta d\zeta
,$$
$$
W_2(t)=\frac{\pi^2x^2}{t^2}\int\limits_{-\infty}^\infty\omega\left(\frac{\zeta}{T_0^{-1+\varepsilon}}\right)
\sinh^2 \zeta\cos(2\pi x\cosh \zeta)\cos 2t\zeta d\zeta
$$ and $B$ arbitrarily large. For simplicity, we only consider the term involving $W_1(t).$ \\
Since $\zeta\leqslant T_0^{-1+\varepsilon},$ by the Taylor expansion of $\cosh \zeta, $
$$
e(x\cosh \zeta)=e(x)e\left(\frac{\zeta^2x}{2}\right)\left(1+\frac{1}{4!}\zeta^4x+\cdots+\frac{1}{(2k)!}\zeta^{2k}x
+O(T_0^{-2k+\varepsilon})\right).
$$
From now on, we only deal with the leading term since all the other terms are similar. We always truncate the series at some point
till the error term is negligible. Now
\begin{equation}
W_1(t)\sim\Im\left\{\frac{\pi x}{2t^2}\int\limits_{-\infty}^\infty\omega\left(\frac{\zeta}{T_0^{-1+\varepsilon}}\right)
\cosh\zeta
e(x)e\left(\frac{\zeta^2x}{2}\right)\cos(2t\zeta)d\zeta\right\}
\end{equation}
which is bounded by $xt^{-2}T_0^{-1+\varepsilon}$ trivially. We are led to estimate
\begin{eqnarray}
\lefteqn{
\sum\limits_{l\geqslant 1}\sum\limits_{m\geqslant 1}\sum\limits_{n\geqslant 1}\frac{A(m, n)}{(m^2nl)^{\frac{1}{2}}}\Omega(l)k(m^2n)
g\left(\frac{l}{N_2}\right)h\left(\frac{m^2n}{N_1}\right)\sqrt{nl}}\\
&&\hspace{1cm}\cdotp \sum\limits_{0<c\leqslant T^{\frac{11}{9}+\varepsilon}m^{-1}}c^{-2}S(n, l; c)
e\left(\frac{2\sqrt{nl}}{c}\right)e\left(\frac{\zeta^2\sqrt{nl}}{c}\right).
\nonumber
\end{eqnarray}
Trivially it is bounded by $N_1N_2T^\varepsilon$ due to Weil's bound (5.6) and (3.12), which yields that
\begin{eqnarray}
\lefteqn{
\sum\limits_{l\geqslant 1}\sum\limits_{m\geqslant 1}\sum\limits_{n\geqslant 1}\frac{A(m, n)}{(m^2nl)^{\frac{1}{2}}}\Omega(l)k(m^2n)
g\left(\frac{l}{N_2}\right)h\left(\frac{m^2n}{N_1}\right)}\\
&&\cdotp \sum\limits_{0<c\leqslant T^{\frac{11}{9}+\varepsilon}m^{-1}}c^{-1}S(n, l; c)
H_1^+\left(\frac{2\sqrt{nl}}{c}\right)\ll N_1N_2T_0^{-1+\varepsilon}.
\nonumber
\end{eqnarray}
In the case that $N_1N_2\leqslant T^{\frac{11}{6}+\varepsilon}T_0,$ the above trivial bound implies that (6.6) is bounded by $T^{\frac{11}{6}+\varepsilon}$, which 
is admissible with the error term in the main theorem. In the following, we assume
\begin{equation}
T_0^{4+\varepsilon}\geqslant N_1N_2\geqslant T^{\frac{11}{6}+\varepsilon}T_0.
\end{equation}
Hence
\begin{equation}
T_0\geqslant T^{\frac{11}{18}+\varepsilon}.
\end{equation}
Depending on the length of $l$, we consider two cases:\\
\begin{equation}
 \hspace{-2cm} 1) \;N_2\leqslant T_0^{\frac{1}{3}};\;\;\;\hspace{3cm} 2)\;N_2\geqslant T_0^{\frac{1}{3}}.
\end{equation}
In this section, we will study the first case. The second case will be done in the next section. Opening the Kloosterman
 sum
$
S(n, l; c)$ as in (2.10) and applying the Voronoi formula Proposition 2.1 for the $n$-sum, we have
\begin{eqnarray}
\lefteqn{\;\;\;\;\;
\sum\limits_{n>0}A(m, n)e\left(\frac{n\bar{d}}{c}\right)
\psi(n)}\\
&&
=\frac{c\pi^{-\frac{5}{2}}}{4i}\sum\limits_{n_1|cm}\sum\limits
_{n_2>0}\frac{A(n_2, n_1)}{n_1n_2}S(md, n_2; m cn_1^{-1})\Psi_{0, 1}^0
\left(\frac{n_2n_1^2}{c^3m}\right)
\nonumber\\
&&\;\;\;\;+\frac{c\pi^{-\frac{5}{2}}}{4i}\sum\limits_{n_1|cm}\sum\limits
_{n_2>0}\frac{A(n_2, n_1)}{n_1n_2}S(m d, -n_2; m cn_1^{-1})
\Psi_{0, 1}^1
\left(\frac{n_2n_1^2}{c^3m}\right)
\nonumber,
\end{eqnarray}
where
\begin{equation}
\Psi_{0,1}^0(x)=\psi_0(x)+\frac{\pi^{-3}c^3m}{n_1^2n_2i}\psi_1(x),\;\; 
\Psi_{0,1}^1(x)=\psi_0(x)-\frac{\pi^{-3}c^3m}{n_1^2n_2i}\psi_1(x),
\end{equation}
for $k=0, 1,$

\begin{equation}
\psi_k(x)=2\pi^3x\int\limits_{(\sigma^{'})}(\pi^3x)^{-2s}G_k(s)\tilde{\psi}(-2s+1-k)ds
\end{equation}
with $\sigma^{'}=\frac{1+\sigma}{2}>\frac{1}{4},$
\begin{equation}
G_k(s)=\frac{
\Gamma\left(s+k+\frac{\alpha}{2}\right)\Gamma\left(s+k+\frac{\beta}{2}\right)\Gamma\left(s+k+\frac{\gamma}{2}\right)}
{
\Gamma\left(\frac{1}{2}-s-\frac{\alpha}{2}\right)
\Gamma\left(\frac{1}{2}-s-\frac{\beta}{2}\right)
\Gamma\left(\frac{1}{2}-s-\frac{\gamma}{2}\right)},
\end{equation}
\begin{equation}
\tilde{\psi}(s)=\int\limits_0^\infty\psi(x)x^s\frac{dx}{x}
\end{equation}
and
\begin{equation}
\psi(x)=e\left(\frac{2\sqrt{xl}}{c}+\frac{\zeta^2\sqrt{xl}}{c}\right)k(m^2x)h\left(\frac{m^2x}{N_1}\right).
\end{equation}
We require an asymptotic formula for $\psi_0(x)$ when $xN_1m^{-2}$ is large. We formulate the asymptotic formula for the general case in the following lemma:
\begin{lemma}
Suppose $\psi(x)$ is a smooth function compactly supported on $[X, 2X],$ $\psi_0(x)$ is defined by (6.12), then for any fixed integer $K\geqslant 1$ and
$xX\gg 1,$ we have
\begin{eqnarray}
\lefteqn{
\psi_0(x)=2\pi^4xi\int\limits_0^\infty\psi(y)\sum\limits_{j=1}^K
\frac{c_j\cos (6\pi x^{\frac{1}{3}}y^{\frac{1}{3}})+d_j\sin (6\pi 
x^{\frac{1}{3}}y^{\frac{1}{3}})}{(\pi^3xy)^{\frac{j}{3}}}dy}\nonumber\\
&&\hspace{6cm}+O\left((xX)^{\frac{-K+2}{3}}\right),
\nonumber
\end{eqnarray}
where $c_j$ and $d_j$ are constants depending on $\alpha, \beta$ and $\gamma,$ in particular, \\$c_1=0, d_1=-\frac{2}{\sqrt{3\pi}}.$
\end{lemma}
\noindent{\bf Proof.} Let 
$$
V(s)=\frac{(-s)\Gamma(3s-1)}{\Gamma\left(\frac{3}{2}-3s\right)}3^{-6s+\frac{5}{2}},
$$
then
$$G_0(s)=V(s)+V(s)H(s)$$
where
$$
H(s)=V(s)^{-1}G_0(s)-1.
$$
Applying Stirling's formula, namely
$$
\log \Gamma(s+c)=\left(s+c-\frac{1}{2}\right)\log s-s +\frac{1}{2}\log 2\pi +\sum\limits_{j=1}^k\frac{a_j}{s^j}+O_\delta\left(\frac{1}{|s|^{K+1}}\right)
$$
which is valid for $c$ a constant, any fixed integer $K\geqslant 1, |\arg s|\leqslant \pi-\delta$ for $\delta>0,$ where
the points $s=0$ and the neighbourhoods of the poles of $\Gamma(s+c)$ are excluded, and the $a_j$ are suitable constants, one shows that
$$
H(s)=\sum\limits_{j=1}^K\frac{b_j}{s^j}+O\left(\frac{1}{|s|^{K+1}}\right)
$$
where $b_j$ are constants depending on $\alpha, \beta$ and $\gamma.$ In the folllowing we will follow closely the proof of Ivic [Iv2]
for the special case $\alpha=\beta=\gamma=0.$ Let
$$
I_1=\frac{2}{2\pi i}\int\limits_{(\sigma^{'})}u^{-2s}V(s)\tilde{\psi}(-2s+1) ds.
$$
Changing variables $3s-1\rightarrow w,$ we have
\begin{eqnarray}
\lefteqn{
I_1=
\frac{-2}{6\pi i}\int\limits_{(\sigma^{''})}(1+w)\frac{\Gamma(w)}{\Gamma(\frac{1}{2}-w)}3^{-\frac{1}{2}-2 w}
u^{-2\frac{w+1}{3}}\tilde{\psi}\left(\frac{-2 w+1}{3}\right) dw}\nonumber\\
&&\hspace{8cm}=I_3+I_4,\nonumber
\end{eqnarray}
where 
$$
I_3=\frac{-2}{6\pi i}\int\limits_{(\sigma^{''})}\frac{\Gamma(w)}{\Gamma(\frac{1}{2}-w)}3^{-\frac{1}{2}-2 w}
u^{-2\frac{w+1}{3}}\tilde{\psi}\left(\frac{-2 w+1}{3}\right) dw
$$
and 
$$
I_4=\frac{-2}{6\pi i}\int\limits_{(\sigma^{''})}w\frac{\Gamma(w)}{\Gamma(\frac{1}{2}-w)}3^{-\frac{1}{2}-2 w}
u^{-2\frac{w+1}{3}}\tilde{\psi}\left(\frac{-2 w+1}{3}\right) dw
$$
with $\sigma^{''}=3\sigma^{'}-1$.
Moving the line of integration in $I_4$ to the left to $\Re s=-\infty,$ we pick up poles of $\Gamma(w)$ at $w=-n$ for
$n=1, 2, \ldots$ with residues $(-1)^n/n!,$ then we have
\begin{eqnarray}
\lefteqn{
I_4=\frac{6}{\sqrt{3}}\sum\limits_{n=1}^\infty\frac{(-1)^n n (3y^{\frac{1}{3}})^{2n-2}\tilde{\psi}\left(\frac{2n+1}{3}\right)}{n! \Gamma\left(n+\frac{1}{2}\right)}
}\nonumber\\
&&\hspace{2cm}=\frac{-2}{\sqrt{3 \pi}}\int\limits_0^\infty\psi(y)\frac{\sin \big(6(uy)^{\frac{1}{3}}\big)}
{(uy)^{\frac{1}{3}}}dy
\nonumber
\end{eqnarray}
%\end{document}
\begin{eqnarray}
\lefteqn{
I_3=\frac{-6}{\sqrt{3}}\int\limits_0^\infty\psi(y)
\frac{J_{-\frac{1}{2}}\big(6(uy)^{\frac{1}{3}}\big)}{\big(3(uy)^{\frac{1}{3}}\big)
^{\frac{3}{2}}}dy}\nonumber\\
&&\hspace{2cm}=
\frac{-2}{3\sqrt{3\pi}}\int\limits_0^\infty\psi(y)\frac{\cos \big(6(uy)^{\frac{1}{3}}\big)}{(uy)^{\frac{2}{3}}}dy
\nonumber\end{eqnarray}
 where we used the well-known integral representation of the $J$-Bessel function (see [EMOT], p.21) 
$$
\frac{1}{2\pi i}\int\limits_{(\sigma)}\frac{\Gamma(s)}{\Gamma(\nu-s+1)}\frac{ds}{x^{2s}}=\frac{J_\nu(2x)}{x^\nu}
$$ for $0<\sigma<\frac{1}{2}\nu+\frac{1}{2}$
and the formula (see [GR], p.914)
$$
J_{-(n+\frac{1}{2})}(x)=\sqrt{\frac{2}{\pi}}x^{n+\frac{1}{2}}\left(\frac{d}{xdx}\right)^n\frac{\cos x}{x}
$$ for $n$ a nonnegative integer and $x>0.$
For $j\geqslant 1,$ let
$$
I_{2, j}=\frac{2}{2\pi i}\int\limits_{(\sigma^{'})}u^{-2s}V(s)\tilde{\psi}(-2s+1) \frac{ds}{s^j}.
$$
Set $3s-1=w,$ then
$$
I_{2, j}=
\frac{-2}{6\pi i}\int\limits_{(\sigma^{''})}(1+w)\frac{\Gamma(w)}{\Gamma(\frac{1}{2}-w)}3^{-\frac{1}{2}-2 w+j}
u^{-2\frac{w+1}{3}}\tilde{\psi}\left(\frac{-2 w+1}{3}\right) \frac{dw}{(w+1)^j}
$$
where $\sigma^{''}=3\sigma^{'}-1.$
%\end{document}
Particularly,
$$
I_{2, 1}=3I_3=\frac{-2}{3\sqrt{3\pi}}\int\limits_0^\infty\psi(y)
\frac{\cos \big(6(uy)^{\frac{1}{3}}\big)}{(uy)^{\frac{2}{3}}}dy
$$
and 
$$
I_{2, 2}=I_{2, 2}^1+I_{2, 2}^2
$$
with
$$
I_{2, 2}^1=\frac{-2}{6\pi i}\int\limits_{(\frac{1}{6})}\frac{-\Gamma(w)}{(\frac{1}{2}-w)
\Gamma(\frac{1}{2}-w)}3^{\frac{3}{2}-2w}
u^{-2\frac{w+1}{3}}\tilde{\psi}\left(\frac{-2w+1}{3}\right)
$$
and
$$
I_{2, 2}^2=\frac{-2}{6\pi i}\int\limits_{(\frac{1}{6})}\frac{3\Gamma(w)}{2(\frac{1}{2}-w)(w+1)
\Gamma(\frac{1}{2}-w)}3^{\frac{3}{2}-2w}
u^{-2\frac{w+1}{3}}\tilde{\psi}\left(\frac{-2w+1}{3}\right)dw.
$$
Using the above integral representation of the $J$-Bessel function and the following formula (see [GR], p.914)
$$
J_{n+\frac{1}{2}}(x)=(-1)^n\sqrt{\frac{2}{\pi}}x^{n+\frac{1}{2}}\left(\frac{d}{xdx}\right)^n\frac{\sin x}{x}
$$ for $n$ a nonnegative integer and $x>0.$
%\end{document}
we have
\begin{eqnarray}
\lefteqn{
I_{2, 2}^1=\frac{2\cdot 3^{\frac{5}{2}}}{9}
\int\limits_0^\infty\psi(y)
\frac{J_{\frac{1}{2}}\big(6(uy)^{\frac{1}{3}}\big)}
{\sqrt{3}(uy)^{\frac{5}{6}}}dy
}\nonumber\\
&&\hspace{1cm}
=\frac{2}{\sqrt{\pi}}\int\limits_0^\infty\psi(y)\frac{\sin\big(6(uy)^{\frac{1}{3}}\big)}{\sqrt{3}(uy)}dy.
\nonumber
\end{eqnarray}
Applying the above precedure repeatly to $I_{2, 2}^2,$ one can derive the lower order terms. The last integral can be estimated
trivially by shifting the line of integration as far as possible. This finishes the proof of the lemma.
\mbox{$\Box$}\\
For later use, we only consider the leading term in Lemma 6.1 since all the other terms can be treated similarly. Now let $\psi(x)$ 
be defined by (6.15).\\
1) In the case $x\frac{N_1}{m^2}\gg T^\varepsilon,$ by the above lemma,
\begin{eqnarray}
\lefteqn{
\psi_0(x)\sim 2\pi^4xi\int\limits_0^\infty\psi(y)\frac{d_1\sin (6\pi x^{\frac{1}{3}}y^{\frac{1}{3}})}
{(\pi^3xy)^{\frac{1}{3}}}dy
}\\
&&=\pi^4xd_1\int\limits_0^\infty
\frac{e\left(\frac{2\sqrt{yl}}{c}+3x^{\frac{1}{3}}y^{\frac{1}{3}}\right)k(m^2y)
h\left(\frac{m^2y}{N_1}\right)
e\left(\frac{\sqrt{yl}\zeta^2}{c}\right)}
{(\pi^3xy)^{\frac{1}{3}}}dy\nonumber\\
&&-\pi^4xd_1\int\limits_0^\infty
\frac{e\left(\frac{2\sqrt{yl}}{c}-3x^{\frac{1}{3}}y^{\frac{1}{3}}\right)k(m^2y)
h\left(\frac{m^2y}{N_1}\right)e\left(\frac{\sqrt{yl}\zeta^2}{c}\right)}
{(\pi^3xy)^{\frac{1}{3}}}dy\nonumber
\end{eqnarray}
%\end{document}

Let 
$$
u_1(y)=\frac{2\sqrt{yl}}{c}+3x^{\frac{1}{3}}y^{\frac{1}{3}},
$$
%\end{document}
then
$$
u_1^{'}(y)=\frac{1}{c}\sqrt{\frac{l}{y}}+x^{\frac{1}{3}}y^{-\frac{2}{3}},
$$
so
$$
u_1^{'}(y)y\geqslant x^{\frac{1}{3}}y^{\frac{1}{3}}\gg T^{\frac{\varepsilon}{3}}.
$$
%\end{document}

By partial integration enough times, one shows that the contribution from the first integral in (6.16) is negligible. Let
$$
u_2(y)=\frac{2\sqrt{yl}}{c}-3x^{\frac{1}{3}}y^{\frac{1}{3}},
$$
then
$$
u_2^{'}(y)=\frac{1}{c}\sqrt{\frac{l}{y}}-x^{\frac{1}{3}}y^{-\frac{2}{3}}.
$$
%\end{document}
1) If $x\geqslant 2\sqrt{\frac{N_1l^3}{m^2c^6}},$ then $|u_2^{'}(y)|\gg y^{-\frac{1}{2}}l^{\frac{1}{2}}c^{-1},$
so \\$|u_2^{'}(y)y|\gg N_1^{\frac{1}{2}}N_2^{\frac{1}{2}}T^{-\frac{11}{9}-\varepsilon}\gg T^\varepsilon$ due to (6.7) and (6.8). In this case,
by partial integration enough times, the contribution from the second integral in (6.16) is negligible. \\
2) If $x\leqslant \frac{2}{3}\sqrt{\frac{N_1l^3}{m^2c^6}},$ then $|u_2^{'}(y)\gg y^{-\frac{1}{2}}l^{\frac{1}{2}}
c^{-1},$ so $$|u_2^{'}(y)y|\gg y^{\frac{1}{2}}l^{\frac{1}{2}}c^{-1}\gg T^\varepsilon.$$
By partial integration enough times, one shows that the contribution from the second integral in (6.16) is negligible.\\
3) If $\frac{2}{3}\sqrt{\frac{N_1l^3}{m^2c^6}}\leqslant x\leqslant 2\sqrt{\frac{N_1l^3}{m^2c^6}},$ then there is 
a stationary phase point $y_0=c^6x^2l^{-3}$ such that $u_2^{'}(y_0)=0.$ Applying the stationary phase method 
([Hu], Lemma 5.5.6), we have
%\end{document}
\begin{eqnarray}
\lefteqn{
\int\limits_{-\infty}^\infty e(u_2(y))e\left(\frac{\zeta^2\sqrt{yl}}{c}\right)k(m^2y)h\left(\frac{m^2y}{N_1}\right)dy}\\
&&=\frac{e(-xc^2l^{-1})e\left(\frac{1}{8}\right)e(\zeta^2c^2xl^{-1})k(m^2c^6x^2l^{-3})h\left(\frac{m^2c^6x^2}{l^3N_1}\right)}{\sqrt{u_2^{''}(y_0)}}
\nonumber\\
&&\hspace{5cm}+ O\left(\frac{c^2T^\varepsilon}{N_2}+\frac{N_1^{\frac{1}{4}}c^{\frac{3}{2}}T^\varepsilon}{m^{\frac{1}{2}}N_2^{\frac{3}{4}}}\right).
\nonumber
\end{eqnarray}
%\end{document}
Since
\begin{equation}
\sum_{\substack{0\leqslant d<c\\ (d, c)=1}}e\left(\frac{ld}{c}\right)S(md, n_2; mcn_1^{-1})=\sum
_{\substack{u
(
\text{mod}\;mcn_1^{-1})\\u\bar{u}\equiv 1 (
\text{mod}\;mcn_1^{-1})}}
S(0, l
+un_1; c)e\left(\frac{n_2\bar{u}}{mcn_1^{-1}}\right)
\end{equation}
where
$$
S(0, a; c)=\sum_{\substack{d(\text{mod}\; c)\\(d, c)=1}}e\left(\frac{ad}{c}\right)
$$
is the Ramanujan sum which is bounded by $(a, c).$ Therefore, (6.18) is bounded by $mc^{1+\varepsilon}.$ 
The contribution to (6.5) from the error term
in (6.17) is bounded by
%\end{document}
\begin{eqnarray}
\lefteqn{\label{(6.19)}}\\
&&\sum\limits_{l\geqslant 1}\sum\limits_{m\geqslant 1}\Omega(l)
g\left(\frac{l}{N_2}\right)
\sum\limits_{0<c\leqslant 
T^{\frac{11}{9}+\varepsilon}m^{-1}}
\sum\limits_{n_1|cm}
\sum\limits_{\frac{2}{3}\frac{\sqrt{N_1l^3}}{n_1^2}\leqslant n_2\leqslant
 2\frac{\sqrt{N_1l^3}}{n_1^2}}
%\frac{|A(n_2, n_1)|}{n_1n_2}
\nonumber
%\end{eqnarray}
%\end{document}
\\
&&\hspace{1cm}.\frac{|A(n_2, n_1)|}{n_1n_2}
\left(\frac{c^2T^\varepsilon}{N_2}+\frac{{N_1}^{\frac{1}{4}}c^{\frac{3}{2}}T^{\varepsilon}}
{m^{\frac{1}{2}}N_2^{\frac{3}{4}}}\right)
\left(\sqrt{\frac{N_1N_2^3}{m^2c^6}}\right)^{\frac{2}{3}}
\left(\frac{N_1}{m^2}\right)^{-\frac{1}{3}}\nonumber\\
&&\hspace{7cm}\ll T_0T^{1+\varepsilon}+T_0^2T^{\frac{1}{2}+\varepsilon}.
\nonumber
\end{eqnarray}
%\end{document}
The contribution to (6.5) from the main term in (6.17) is bounded by 
\begin{eqnarray}
\lefteqn{\;\;\;\;
\sum\limits_{l\geqslant 1}\sum\limits_{m\geqslant 1}\Omega(l)g\left(\frac{l}{N_2}\right)
\!\!\!\sum\limits_{0<c\leqslant T^{\frac{11}{9}+\varepsilon}m^{-1}}
\sum\limits_{n_1|cm}
\!\sum\limits_{\frac{2}{3}\frac{\sqrt{N_1l^3}}{n_1^2}\leqslant n_2\leqslant 2
\frac{\sqrt{N_1l^3}}{n_1^2}}
}\\
%\sum\limits_{\frac{1}{1000}\frac{\sqrt{N_1N_2^3}}{n_1^2}\leqslant n_2\leqslant 2
%\frac{\sqrt{N_1N_2^3}}{n_1^2}}}\\
&&\hspace{1.5cm}
%\sum\limits_{n_1|cm}
%\sum\limits_{\frac{1}{1000}\frac{\sqrt{N_1N_2^3}}{n_1^2}\leqslant n_2\leqslant 1000
%\frac{\sqrt{N_1N_2^3}}{n_1^2}}
\cdotp\frac{|A(n_1, n_2)|}{n_1n_2}
l^{-\frac{5}{2}}c^5\left(\frac{N_1N_2^3}{m^2c^6}\right)^{\frac{3}{4}+\frac{1}{3}}
\left(\frac{N_1}{m^2}\right)^{-\frac{1}{3}}
\ll T_0^{\frac{17}{6}+\varepsilon}
\nonumber
%&&\hspace{5cm}\ll T_0^{\frac{17}{6}+\varepsilon}
\end{eqnarray}
where we used the condition that
\begin {equation}
N_1\ll T_0^{3+\varepsilon},\;\;\;\;\;\;\;\;\;
N_2\ll T_0^{\frac{1}{3}}
.\end{equation}
Now if $x\frac{N_1}{m^2}\ll T^{\varepsilon}$ (recall $\psi_0(x)$ is defined by (6.12)), let $\sigma^{'}=\frac{1}{3}$
\begin{eqnarray}
\lefteqn{
\psi_0(x)\ll x\int\limits_{(\sigma^{'})}x^{-2\sigma^{'}}(|s|+1)^{6\sigma^{'}-\frac{3}{2}}\frac{\left(\frac{N_1}{m^2}\right)^{-2\sigma^{'}+1}}
{(|s|+1)^{100}}ds}\nonumber\\
&&\hspace{3cm}\ll x^{1-2\sigma^{'}}\left(\frac{N_1}{m^2}\right)^{-2\sigma^{'}+1}\ll T^{\frac{\varepsilon}{3}}.
\nonumber
\end{eqnarray}
%\end{document}
whose contribution to (6.5) is bounded by
\begin{eqnarray}
\lefteqn{
\sum\limits_{l\geqslant 1}\sum\limits_{m\geqslant 1}\sum\limits_{0<c\leqslant T^{\frac{11}{9}+\varepsilon}m^{-1}}\sum\limits_{n_1|cm}
\sum\limits_{n_2\leqslant \frac{T^\varepsilon m^3c^3}{n_1^2N_1}}\frac
{|A(n_2, n_1)|}{n_1n_2}g\left(\frac{l}{N_2}\right)T^{\frac{\varepsilon}{3}}}\\
&&\hspace{6cm}
\ll N_2T^{\frac{11}{9}+\varepsilon},
\nonumber
\end{eqnarray}
where we used (3.13) again.
Combining (6.19), (6.20) and (6.22), we conclude that under the condition (6.21), (6.6) is bounded by \\
$(T_0T^{1+\varepsilon}+T_0^2T^{\frac{1}{2}+\varepsilon}+T_0^2T^{\frac{1}{2}+\varepsilon}+T_0^{\frac{17}{6}+\varepsilon})T_0^{-1}
\ll T^{\frac{11}{6}+\varepsilon}$ which is admissible with the error term in the main theorem.
\section{Sums of Kloosterman sums - small $c$: Part II}
\setcounter{equation}{0} 
In this section, we continue to estimate the contribution from sums of Kloosterman sums for small $c$, i.e., (6.5), under the condition
that 
\begin{equation}
N_1\leqslant T^{3+\varepsilon},\;\; N_1N_2\geqslant T^{\frac{11}{6}+\varepsilon},\;\; T_0^{\frac{1}{3}}\leqslant N_2\leqslant T_0^{1+\varepsilon}.
\end{equation}
Opening the Kloosterman sum as in (2.10) and applying the Poisson summation for the $l$-sum in (6.5), we have 
\begin{eqnarray}
\lefteqn{\;\;\;\;
\sum\limits_{l\in \mathbb{Z}}e\left(\frac{ld+2\sqrt{nl}}{c}\right)e\left(\frac{\zeta^2\sqrt{nl}}{c}\right)\Omega(l)g\left(\frac{l}{N_2}\right)
}\\
&&=\sum\limits_{k\in\mathbb{Z}}\int\limits_{-\infty}^{\infty}e\left(\frac{(kc+d)x+2\sqrt{nx}}{c}\right)e\left(\frac{\zeta^2\sqrt{nx}}{c}\right)\Omega(x)g\left(\frac{x}{N_2}\right)dx.
\nonumber
\end{eqnarray}
Let
$$
w(x)=\frac{(kc+d)x+2\sqrt{nx}}{c},
$$
then
$$
w{'}(x)=\frac{(kc+d)+\sqrt{\frac{n}{x}}}{c}.
$$
There are two cases:\\
1) For $|kc+d|\geqslant 10\sqrt{\frac{N_1}{N_2m^2}},$ then
$$
w^{'}(x)N_2\gg\frac{1}{c}\sqrt{\frac{N_1N_2}{m^2}}\gg T^{\varepsilon}.
$$
By partial integration $\left[\frac{A}{\varepsilon}\right]+1$ times
\begin{eqnarray}
\lefteqn{
\int\limits_0^\infty e(w(x))e\left(\frac{\zeta^2\sqrt{nx}}{c}\right)\Omega(x)g\left(\frac{x}{N_2}\right)dx}\nonumber\\
&&\hspace{2cm}
\ll \frac{N_2}{\left(\frac{|kc+d|N_2}{c}\right)^{\frac{A}{\varepsilon}}}\ll N_2T^{-A}
\nonumber
\end{eqnarray}
where $A>0$ is arbitrarily large. Thus the contribution from such terms is negligible.\\
2) For $|kc+d|\leqslant 
\frac{1}{10}\sqrt{\frac{N_1}{N_2m^2}},$
$$
w^{'}(x)N_2\gg \frac{\sqrt{N_1N_2m^{-2}}}{c}\gg T^{\varepsilon},
$$
as the above, by partial integration $\left[\frac{A}{\varepsilon}\right]+1$ times,
\begin{eqnarray}
\lefteqn{
\int\limits_0^\infty e(w(x)e\left(\frac{\zeta^2\sqrt{nx}}{c}\right)\Omega(x)g\left(\frac{x}{N_2}\right)dx}\nonumber\\
&&\hspace{2cm}
\ll \frac{N_2}{\left(\frac{1}{c}\sqrt{\frac{N_1N_2}{m^2}}\right)^{\frac{A}{\varepsilon}}}\ll N_2T^{-A}
\nonumber
\end{eqnarray}
where $A>0$ is arbitrarily large. Hence the contribution from such terms is negligible.\\
3) For $\frac{1}{10}\sqrt{\frac{N_1}{N_2m^2}}\leqslant
|kc+d|\leqslant 10\sqrt{\frac{N_1}{N_2m^2}},$
there is a stationary phase point $x_0=\frac{n}{(kc+d)^2}$ such that $w^{'}(x_0)=0.$ By the stationary phase method 
(see [Hu], Lemma 5.5.6), we have
\begin{eqnarray}
\lefteqn{
\int\limits_0^\infty e(w(x)e\left(\frac{\zeta^2\sqrt{nx}}{c}\right)\Omega(x)g\left(\frac{x}{N_2}\right)dx=\sqrt{2cn}e\left(\frac{-n}{c^2k+cd}+\frac{1}{8}\right)
}\nonumber\\
&&
%=\sqrt{2cn}e\left(\frac{-n}{c^2k+cd}+\frac{1}{8}\right)
\hspace{1cm}\cdotp e\left(\frac{-\zeta^2n}{c^2k+cd}\right)|kc+d|^{-\frac{3}{2}}
\Omega\left(\frac{n}{(kc+d)^2}\right)g\left(\frac{n}{(kc+d)^2N_2}\right)\nonumber\\
&&\hspace{4cm}+O_{\varepsilon}\left(c^2m^2T^{\varepsilon}N_1^{-1}+c^{\frac{3}{2}}N_2^{\frac{1}{4}}N_1^{-\frac{3}{4}}m^{\frac{3}{2}}T^{\varepsilon}\right).
\nonumber
\end{eqnarray}
(7.2) becomes
\begin{eqnarray}
\lefteqn{
\sum\limits_{\frac{1}{10}\sqrt{\frac{N_1}{N_2m^2}}\leqslant |kc+d|\leqslant 10\sqrt{\frac{N_1}{N_2m^2}}}\sqrt{2cn}|kc+d|^{-\frac{3}{2}}
e\left(\frac{-\zeta^2n}{c^2k+cd}\right)
}\\
&&\hspace{1cm}\cdotp e\left(\frac{-n}{c^2k+cd}+\frac{1}{8}\right)\Omega\left(\frac{n}{(kc+d)^2}\right)g\left(
\frac{n}{(kc+d)^2N_2}\right)\nonumber\\
&&\hspace{1cm}+O_{\varepsilon}\left(\frac{cmT^\varepsilon}{\sqrt{N_1N_2}}+\frac{c^{\frac{1}{2}}m^{\frac{1}{2}}T^{\varepsilon}}
{N_2^{\frac{1}{4}}N_1^{\frac{1}{4}}}+c^2m^2T^{\varepsilon}N_1^{-1}+c^{\frac{3}{2}}N_2^{\frac{1}{4}}
N_1^{-\frac{3}{4}}m^{\frac{3}{2}}
T^{\varepsilon}\right)
.
\nonumber
\end{eqnarray}
The contribution to (6.5) from the above error term is bounded by 
\begin{eqnarray}
\lefteqn{
\sum\limits_m\sum\limits_n\frac{|A(m, n)|}{m}\sum\limits_{0<c\leqslant T^{\frac{11}{9}+\varepsilon}m^{-1}}c^{-1}}\\
&&\cdotp\left[\frac{cmT^\varepsilon}{\sqrt{N_1N_2}}+\frac{c^{\frac{1}{2}}m^{\frac{1}{2}}T^\varepsilon}{N_2^{\frac{1}{4}}
N_1^{\frac{1}{4}}}
+c^2m^2T^\varepsilon N_1^{-1}+c^{\frac{3}{2}}N_2^{\frac{1}{4}}N_1^{-\frac{3}{4}}m^{\frac{3}{2}}T^{\varepsilon}\right]\nonumber\\
&&
\leqslant 
N_1^{\frac{1}{2}}N_2^{-\frac{1}{2}}T^{\frac{11}{9}+\varepsilon}+N_1^{\frac{3}{4}}N_2^{-\frac{1}{4}}T^{\frac{11}{18}+\varepsilon}
+T^{\frac{22}{9}+\varepsilon}+N_1^{\frac{1}{4}}N_2^{\frac{1}{4}}T^{\frac{11}{6}+\varepsilon}.
\nonumber
\end{eqnarray}
%\end{document}
Therefore the contribution to (6.6) from the above error term is:
\begin{eqnarray}
\lefteqn{\hspace{-1cm}\;\;\;\;
O\big((N_1^{\frac{1}{2}}N_2^{-\frac{1}{2}}T^{\frac{11}{9}+\varepsilon}+N_1^{\frac{3}{4}}N_2^{-\frac{1}{4}}T^{\frac{11}{18}+\varepsilon}
+T^{\frac{22}{9}+\varepsilon}+N_1^{\frac{1}{4}}N_2^{\frac{1}{4}}T^{\frac{11}{6}+\varepsilon})T_0^{-1+\varepsilon}\big)}
\nonumber\\
&&\hspace{7cm}=O_{\varepsilon, f}
\left(T^{\frac{11}{6}+\varepsilon}\right).\nonumber
\end{eqnarray}
The contribution to (6.5) from the main term in (7.3) is
\begin{eqnarray}
\lefteqn{\;\;\;\;\;\;\;\;
\sum\limits_{m\geqslant 1}\sum\limits_{n\geqslant 1}\frac{A(m, n)}{m}k(m^2n)h\left(\frac{m^2n}{N_1}\right)\sum\limits_
{1\leqslant c\leqslant T^{\frac{11}{9}+\varepsilon}m^{-1}}c^{-2}\sum_{\substack{0\leqslant d<c\\(d, c)=1}}}\\
&&\cdotp\;\sum\limits_{\frac{1}{10}\sqrt{\frac{N_1}{N_2m^2}}\leqslant |kc+d|\leqslant 10\sqrt{\frac{N_1}{N_2m^2}}}\!\!\!\!\!
\sqrt{2cn}|kc+d|^{-\frac{3}{2}}e\left(\frac{-\zeta^2n}{c^2k+cd}+\frac{1}{8}\right)\nonumber\\
&&\hspace{2cm}\cdotp\Omega\left(\frac{n}{(kc+d)^2}\right)
g\left(\frac{n}{(kc+d)^2N_2}\right)e\left(\frac{n\bar{d}}{c}+\frac{-n}{c^2k+cd}\right).
\nonumber
\end{eqnarray}
For the $n$-sum, we apply the Voronoi formula  on ${\rm GL}(3),$ i.e. Proposition 2.1:
\begin{equation}
\sum\limits_{n\geqslant 1}A(m, n)e\left(\frac{n\bar{d}}{c}-\frac{n}{c^2k+cd}\right)\phi(n)
\end{equation}
where
$$
\phi(x)=\sqrt{x}e\left(\frac{-\zeta^2x}{c^2k+cd}\right)k(m^2x)h\left(\frac{m^2x}{N_1}\right)\Omega\left(\frac{x}{(kc+d)^2}\right)
g\left(\frac{x}{((kc+d)^2N_2}\right).
$$
Since
$$
\frac{\bar{d}}{c}-\frac{1}{c^2k+cd}=\frac{\bar{d}(ck+d)-1}{c(ck+d)}:=\frac{d^{'}}{c^{'}},
$$
obviously $c^{'}|ck+d,$ by the Voronoi formula, (7.6) is equal to
\begin{eqnarray}
\lefteqn{\;\;\;\;\;\;
\frac{c^{'}\pi^{-\frac{5}{2}}}{4i}\sum\limits_{n_1|c^{'}m}\sum\limits_{n_2>0}\frac{A(n_2, n_1)}{n_1n_2}
S(m\bar{d^{'}}, n_2; mc^{'}n_1^{-1})
\Phi_{0, 1}^0\left(\frac{n_2n_1^2}{{c^{'}}^3m}\right)}\\
&&+\frac{c^{'}\pi^{-\frac{5}{2}}}{4i}\sum\limits_{n_1|c^{'}m}\sum\limits_{n_2>0}\frac{A(n_2, n_1)}{n_1n_2}
S(m\bar{d^{'}}, -n_2; mc^{'}n_1^{-1})
\Phi_{0, 1}^1\left(\frac{n_2n_1^2}{{c^{'}}^3m}\right),
\nonumber
\end{eqnarray}
where $\Phi_{0, 1}^0(x)$ and $\Phi_{0, 1}^1(x)$ are defined by (2.13) and (2.14), respectively.\\
We only consider the contribution from $\Phi_0(x)$ (recall (2.12)) since all the other terms can be estimated
in the same way.\\
By making a change of the variable,
$$
\Phi_0(x)=2\pi^3x\int\limits_{(\sigma^{'})}(\pi^3x)^{-2s}G(s)\tilde{\phi}(-2s+1)ds
$$
where $G(s)$ is defined by (6.13) and $\sigma^{'}=\frac{1+\sigma}{2}>\frac{1}{4}.$ As in the last section, we consider two cases
seperately.\\
1) When $x\frac{N_1}{m^2}\geqslant T^\varepsilon:$\\
Lemma 6.1 yields that
$$
\Phi_0(x)\sim 2\pi^4ixd_1\int\limits_0^\infty\phi(y)\sin(6\pi x^{\frac{1}{3}}y^{\frac{1}{3}})(\pi^3xy)^{-\frac{1}{3}}dy.
$$
By partial integration $\left[\frac{A}{\varepsilon}\right]+1$ times with $A$ arbitrarily large, we obtain
$$
\Phi_0(x)\ll x^{\frac{2}{3}}\left(\frac{N_1}{m^2}\right)^{\frac{7}{6}}T^{-\frac{A}{3}}
$$
whose contribution to (7.6) is negligible.\\
2) When $x\frac{N_1}{m^2}\leqslant T^\varepsilon:$ we set $\sigma^{'}=\frac{1}{3}$\\
Since
$$
G(s)\ll (|s|+1)^{6\sigma^{'}-\frac{3}{2}}
$$
and 
$$
\tilde{\phi}(-2s+1)\ll \frac{\left(\frac{N_1}{m^2}\right)^{\frac{3}{2}-2\sigma^{'}}}{(|s|+1)^{100}},
$$
we have
\begin{eqnarray}
\lefteqn{
\Phi_0(x)\ll x\int\limits_{(\sigma^{'})}x^{-2\sigma^{'}}(|s|+1)^{6\sigma^{'}-\frac{3}{2}}\frac{\left(\frac{N_1}{m^2}\right)^{\frac{3}{2}-2\sigma^{'}}}
{(|s|+1)^{100}}ds}
%\ll x^{1-2\sigma^{'}}\left(\frac{N_1}{m^2}\right)^{\frac{3}{2}-2\sigma^{'}}\leqslant 
%T^{\frac{\varepsilon}{3}}\left(\frac{N_1}{m^2}\right)^{\frac{1}{2}}
\nonumber\\
&&\hspace{3cm}
\ll x^{1-2\sigma^{'}}\left(\frac{N_1}{m^2}\right)^{\frac{3}{2}-2\sigma^{'}}\leqslant T^{\frac{\varepsilon}{3}}\left(\frac{N_1}{m^2}\right)^{\frac{1}{2}}
\nonumber
\end{eqnarray}
In this case, by Weil's bound (5.6) and (7.7), (7.6) is bounded by 
\begin{eqnarray}
\lefteqn{
c^{'}\sum\limits_{n_1|c^{'}m}\sum\limits_{n_2\leqslant \frac{{c^{'}}^3m^3T^\varepsilon}{n_1^2N_1}}\frac{|A(n_2, n_1)|}{n_1n_2}
(m\bar{d^{'}}, n_2, mc^{'}n_1^{-1})^{\frac{1}{2}}}\nonumber\\
&&\hspace{1cm}\cdotp (mc^{'}n_1^{-1})^{\frac{1}{2}+\varepsilon}\left(\frac{n_2n_1^2}{{c^{'}}^3m}\right)^{1-2\sigma^{'}}\left(\frac{N_1}{m^2}\right)^{\frac{3}{2}-2\sigma^{'}}
\ll {c^{'}}^{\frac{3}{2}}N_1^{\frac{1}{2}}
\nonumber
%&&\hspace{4cm}
%\ll {c^{'}}^{\frac{3}{2}}N_1^{\frac{1}{2}}
%\nonumber
\end{eqnarray}
by (3.13) and the partial summation formula.\\
(7.5) is bounded by 
\begin{eqnarray}
\lefteqn{\hspace{-2cm}
\sum\limits_{m\geqslant 1}\frac{1}{m}\sum\limits_{1\leqslant c\leqslant T^{\frac{11}{9}+\varepsilon}m^{-1}}
c^{-\frac{1}{2}}\!\!\!\!\!\sum\limits_
{\frac{1}{10}\sqrt{\frac{N_1}{N_2m^2}}\leqslant |kc+d|\leqslant 10\sqrt{\frac{N_1}{N_2m^2}}}|kc+d|^{-\frac{3}{2}}
{c^{'}}^{\frac{3}{2}}N_1^{\frac{1}{2}}}\nonumber\\
&&\hspace{3cm}\ll N_1N_2^{-\frac{1}{2}}+N_1^{\frac{1}{2}}T^{\frac{11}{18}+\varepsilon}
\nonumber
\end{eqnarray}
whose contribution to (6.6) is
$$
O\big((N_1N_2^{-\frac{1}{2}}+N_1^{\frac{1}{2}}T^{\frac{11}{18}+\varepsilon})T_0^{-1+\varepsilon}\big)=O(T^{\frac{11}{6}+\varepsilon})
,$$
where we used the condition (7.1). The above error term is admissible with the error term in the main theorem.
\section{Bilinear forms of Kloosterman sums}
\setcounter{equation}{0} 
In this section, we will study the contribution from sums of Kloosterman sums for $c$ in the intermediate range: $T^{\frac{11}{9}+\varepsilon}m^{-1}\leqslant c\leqslant N_1^{\frac{1}{2}}N_2^{\frac{1}{2}}m^{-1}.$ We
split the $m$-sum into two ranges:\\
I) $m\geqslant \sqrt{\frac{N_1}{N_2}}T^{-\varepsilon};$\; \;\;\;\;\;\;\; \hspace{3cm}II) $m\leqslant\sqrt{\frac{N_1}{N_2}}T^{-\varepsilon}.$\\
$\bullet$ For the first range, we apply Weil's bound (5.6) and (3.13),\\
\begin{eqnarray}
\lefteqn{
\sum\limits_{l\geqslant 1}\sum\limits_{m\geqslant \sqrt{\frac{N_1}{N_2}}T^{-\varepsilon}}\sum\limits_{n\geqslant 1}
\frac{A(m, n)}{m}\Omega(l)k(m^2n)g\left(\frac{l}{N_2}\right)h\left(\frac{m^2n}{N_1}\right)}\\
&&\cdotp\sum\limits_{T^{\frac{11}{9}+\varepsilon}m^{-1}\leqslant c\leqslant N_1^{\frac{1}{2}}N_2^{\frac{1}{2}}m^{-1}}
c^{-2}S(n, l; c)e\left(\frac{2\sqrt{nl}}{c}\right)e\left(\frac{\zeta^2\sqrt{nl}}{c}\right)\nonumber\\
&&\hspace{5cm}\ll N_1^{\frac{3}{4}}N_2^{\frac{5}{4}}T^{-\frac{11}{18}}\ll T_0^{\frac{7}{2}}T^{-\frac{11}{18}}\nonumber
\end{eqnarray}
$\bullet$ For the second range, we consider the following bilinear forms which techniques are used in [DI]:
\begin{equation}
\sum\limits_{n\geqslant 1}\sum\limits_{l\geqslant 1}a(n) b(l) S(n, l; c)e\left(\frac{\theta\sqrt{nl}}{c}\right)h\left(\frac{m^2n}{N}\right)
k(m^2n)
\end{equation}
where
$$
a(n)=A(m, n)k(m^2n),\;\;\;b(l)=\Omega(l)g\left(\frac{l}{N_2}\right),\;\;\;\theta=2+\zeta^2.
$$
By Cauchy's inequality, (8.2) is bounded by 
\begin{eqnarray}
\lefteqn{\label{(8.3)}}\\
%\left[\sum\limits_{n\geqslant 1}|a(n)|^2h\left(\frac{m^2n}{N_1}\right)\right]^{\frac{1}{2}}
&&\!\!\!\left[\sum\limits_{n\geqslant 1}|a(n)|^2h\left(\frac{m^2n}{N_1}\right)\right]^{\frac{1}{2}}
\left[\sum\limits_{n\geqslant 1}\left|\sum\limits_{l\geqslant 1}b(l)S(n, l; c)e\left(\frac{\theta\sqrt{nl}}{c}\right)\right|^2h\left(\frac{m^2n}{N_1}\right)\right]^{\frac{1}{2}}
\nonumber\\
&&\leqslant
 %N_1^{\frac{1}{2}}
\left[\sum\limits_{l_1\geqslant 1}\sum\limits_{l_2\geqslant 1}b(l_1)\bar{b(l_2)}
\sum_{\substack{0\leqslant d_1<c\\(d_1, c)=1}}\sum_{\substack{0\leqslant d_2<c\\(d_2, c)=1}}e\left(\frac{l_1\bar{d_1}-l_2\bar{d_2}}{c}\right)
\sum\limits_{n\in\mathbb{Z}}
F(n)
%e(\eta(n))h\left(\frac{m^2n}{N_1}\right)
\right]^{\frac{1}{2}}
N_1^{\frac{1}{2}}
\nonumber
\end{eqnarray}
where
$$ F(n)=e(p(n))h\left(\frac{m^2n}{N_1}\right)$$
with
$$
p(x)=\frac{d_1-d_2}{c}x+\frac{2(\sqrt{l_1}-\sqrt{l_2})\sqrt{x}}{c}\theta.
$$
By the Poisson summation formula,
\begin{equation}
\sum\limits_{n\in \mathbb{Z}}e(p(n))h\left(\frac{m^2n}{N_1}\right)=\sum\limits_{k\in\mathbb{Z}}\int\limits_{\mathbb{R}}
e(p(x)-kx)h\left(\frac{m^2x}{N_1}\right)dx.
\end{equation}
Let $A=\frac{d_1-d_2}{c},$ $B=\frac{2(\sqrt{l_1}-\sqrt{l_2})}{c}\theta.$ If $k\neq A,$ then $|k-A|\geqslant \frac{1}{c}.$
Since
$$
\frac{\sqrt{l_1}-\sqrt{l_2}}{c\sqrt{n}}\theta\ll\frac{l_1-l_2}{c\sqrt{\frac{N_1N_2}{m^2}}}\ll \frac{m}{c}\sqrt{\frac{N_2}{N_1}}
\ll\frac{T^{-\varepsilon}}{c},
$$
$$
|k-p^{'}(x)|\frac{N_1}{m^2}\gg |k-A|\frac{N_1}{m^2}\gg\frac{1}{c}\frac{N_1}{m^2}\gg T^{\varepsilon}.
$$
By partial integration $p$-times,
$$
\int\limits_{\mathbb{R}}e(p(x)-kx)h\left(\frac{m^2x}{N_1}\right)dx\ll\frac{N_1}{m^2}\left(|A-k|\frac{N_1}{m^2}\right)^{-p}.
$$
Hence on taking $p=[B/\varepsilon]+1$ with $B$ arbitrarily large, we deduce that
\begin{equation}
\sum\limits_{k\neq A}\int\limits_{\mathbb{\mathbb{R}}}e(p(x)-kx)h\left(\frac{m^2x}{N_1}\right)dx\ll \frac{N_1}{m^2}T^{-B}
\end{equation}
which is negligible.\\
If $k=A,$ then $k=A=0,$ $d_1=d_2.$
\begin{equation}
\int\limits_{\mathbb{R}}e(p(x))h\left(\frac{m^2x}{N_1}\right)dx\ll
\left\{
\begin{array}{lll}
\frac{N_1}{m^2}& \mbox{} &\text{if} \;\;\;l_1=l_2 \\
\frac{c\sqrt{\frac{N_1N_2}{m^2}}}{|l_1-l_2|}&\mbox{}&  \text{if}\;\;\;l_1\neq l_2 \end{array} \right.
\end{equation}
where we used partial integration once in the case that $l_1\neq l_2.$ Combining (8.3), (8.5) and (8.6), it yields that
(8.2) is bounded by
%\end{document}
\begin{eqnarray}
\lefteqn{\hspace{-2cm}
%N_1^{\frac{1}{2}}
\left[\sum\limits_{l\geqslant 1}|b(l)|^2c\frac{N_1}{m^2}+\sum\limits_{l\geqslant 1}
\sum\limits_{1\leqslant l_2\neq l_1}|b(l_1)||b(l_2)||S(0, l_1-l_2; c)|
\frac{c\sqrt{\frac{N_1N_2}{m^2}}}{|l_1-l_2|}\right]^{\frac{1}{2}}N_1^{\frac{1}{2}}
}\nonumber\\
&&\hspace{6cm} \leqslant \frac{N_1}{m}c^{\frac{1}{2}}N_2^{\frac{1}{2}},
\nonumber
\end{eqnarray}
where we used (3.12).
%\end{document}
It yields that
\begin{eqnarray}
\lefteqn{
\sum\limits_{l\geqslant 1}\sum\limits_{m\leqslant \sqrt{\frac{N_1}{N_2}}T^{-\varepsilon}}\sum\limits_{n\geqslant 1}\frac{A(m, n)}
{m}\Omega(l)k(m^2n)g\left(\frac{l}{N_2}\right)h\left(\frac{m^2n}{N_1}\right)}\\
&& \cdotp\sum\limits_{T^{\frac{11}{9}+\varepsilon}m^{-1}\leqslant c\leqslant N_1^{\frac{1}{2}}N_2^{\frac{1}{2}}m^{-1}}
c^{-2}S(n, l; c)e\left(\frac{2\sqrt{nl}}{c}\right)e\left(\frac{\zeta^2\sqrt{nl}}{c}\right)
\nonumber\\
&&\hspace{5cm}\ll N_1N_2^{\frac{1}{2}}T^{-\frac{11}{18}}\ll T_0^{\frac{7}{2}}T^{-\frac{11}{18}}.
\nonumber
\end{eqnarray}
Gathering (8.1) and (8.7), we conclude that
\begin{eqnarray}
\lefteqn{
\sum\limits_{l\geqslant 1}\sum\limits_{m\geqslant 1}\sum\limits_{n\geqslant 1}\frac{A(m, n)}{m}\Omega(l)k(m^2n)g\left(\frac{l}{N_2}\right)
h\left(\frac{m^2n}{N_1}\right)}\nonumber\\
&& \cdotp\sum\limits_{T^{\frac{11}{9}+\varepsilon}m^{-1}\leqslant c\leqslant N_1^{\frac{1}{2}}N_2^{\frac{1}{2}}m^{-1}}c^{-1}S(n, l; c)
H_1^{+}\left(\frac{2\sqrt{nl}}{c}\right)\nonumber\\
&& \hspace{5cm}\ll T_0^{\frac{5}{2}}T^{-\frac{11}{18}}\ll T^{\frac{11}{9}+\varepsilon}\nonumber
\end{eqnarray}
which is admissible with the error term in the main theorem. This finishes the proof of the main theorem. \mbox{$\Box$}

\;
\;
\;
\noindent {\bf Acknowledgements}\\

\noindent The author would like to thank Dorian Goldfeld for many illuminating discussions and for his encouragements. She would also 
like to thank Peter Sarnak, Wenzhi Luo, Akshay Venkatesh, Jianya Liu and the referee for valuable comments.\\
%\end{document}

\;
\;
%\noindent Dorian Goldfeld: Department of Mathematics, Columbia University, 2990 Broadway, New York, NY, 10027\\
%E-mail address: dg@math.columbia.edu\\

%\noindent Xiaoqing Li: Department of Mathematics, Columbia University, 2990 Broadway, New York, NY, 10027\\
%E-mail address: xqli@math.columbia.edu         

\end{document}